\documentclass{article}
\PassOptionsToPackage{hyphens}{url}
\usepackage[lang=american]{ems-journal}
\usepackage{mathtools,microtype}
\usepackage{tikz}
\usepackage{placeins}
\numberwithin{equation}{section}
\theoremstyle{plain}
\newtheorem{theorem}{Theorem}[section]
\newtheorem{proposition}[theorem]{Proposition}
\newtheorem{lemma}[theorem]{Lemma}
\newtheorem{corollary}[theorem]{Corollary}
\newtheorem{smoothinput}{Theorem}

\theoremstyle{definition}
\newtheorem{definition}[theorem]{Definition}
\theoremstyle{remark}
\newtheorem{remark}[theorem]{Remark}
\newcommand{\R}{\mathbb R}
\newcommand{\D}{\mathbb D}
\newcommand{\Sone}{\mathbb S^1}

\DeclareMathOperator{\dist}{dist}
\hypersetup{pdftitle={Hot Spots on Simply Connected Lipschitz Domains},pdfauthor={Rodrigo Banuelos}}
\begin{document}
\title{Hot Spots on Simply Connected Lipschitz Domains}
\titlemark{Hot Spots on Simply Connected Lipschitz Domains}
\emsauthor{1}{\givenname{Rodrigo}\surname{Ba\~nuelos}}{R.~Ba\~nuelos}
\Emsaffil{1}{
  \department{Department of Mathematics}
  \organisation{Purdue University}
  \zip{47907}
  \city{West Lafayette, IN}
  \country{USA}
  \affemail{banuelos@purdue.edu}}

\begin{abstract}
Using the smooth-domain theorem of OpenAI (see also \cite{BanuelosSmooth2026}) we prove strict boundary
attainment of the extrema of every first positive Neumann eigenfunction
on a bounded simply connected planar Lipschitz domain. No simplicity
assumption is imposed. The proof selects a prescribed eigendirection by
a thin-strip perturbation and then uses conformal exhaustion and a
critical-set argument. We also give polygonal simplifications and
explicit nonconvex starlike examples.
\end{abstract}
\maketitle
\thispagestyle{empty}
%

\section{Introduction and statement of main results}
\label{sec:introduction}
Let $\Omega\subset\R^2$ be a bounded Lipschitz domain. The definition is
recalled in Section~\ref{subsec:conventions}. Let $-\Delta_\Omega$ be the
nonnegative self-adjoint Neumann operator associated with the form
$v\mapsto\int_\Omega|\nabla v|^2$ on $H^1(\Omega)$.
We denote its first positive eigenvalue by $\Lambda_\Omega^{\mathrm N}$
and the corresponding eigenspace by
\[
 E_1(\Omega)=\ker(-\Delta_\Omega-\Lambda_\Omega^{\mathrm N}).
\]
Every member of $E_1(\Omega)$ extends continuously to
$\overline\Omega$ \cite[Theorem~3.14(ii)]{Nittka2011}.

The hot spots conjecture goes back to a lecture given by J.~Rauch at
Tulane University in 1974. The heat-conduction problem discussed there
is in \cite{Rauch1975}. The conjecture asserts that the functions in
$E_1(\Omega)$ attain their maximum and their minimum over
$\overline\Omega$ on $\partial\Omega$. It has a strict form and a
non-strict form. The formulations and the early history are discussed
in \cite[Section~1]{BB}. Following that discussion, we say
that $\Omega$ has \emph{property (HS2)} if every $u\in E_1(\Omega)$
satisfies
\begin{equation}\label{eq:main-HS2}
 \max_{\overline\Omega}u=\max_{\partial\Omega}u,
 \qquad
 \min_{\overline\Omega}u=\min_{\partial\Omega}u.
\end{equation}
We say that $\Omega$ has \emph{property (HS1)} if
\begin{equation}\tag{HS1}\label{eq:HS1}
 \min_{\partial\Omega}u<u(x)<\max_{\partial\Omega}u
 \qquad(x\in\Omega,\ 0\ne u\in E_1(\Omega)).
\end{equation}
Property (HS2) allows an extreme value to be attained both on the
boundary and in the interior. Property (HS1) does not. The strict form
of the conjecture for $\Omega$ is property (HS1). The non-strict form
is property (HS2).

\begin{theorem}[Lipschitz domains]
\label{thm:lipschitz}
Every bounded simply connected Lipschitz domain $\Omega\subset\R^2$ has
property (HS1). That is, every nonzero $u\in E_1(\Omega)$ satisfies
\eqref{eq:HS1}.
\end{theorem}

No assumption is made on the multiplicity of
$\Lambda_\Omega^{\mathrm N}$.
Theorem~\ref{thm:lipschitz} applies to simply connected polygonal domains,
bounded simply connected $C^{1,\alpha}$ domains, $0<\alpha<1$, and
bounded starlike Lipschitz domains, including nonconvex examples.

The smooth input is the following theorem of OpenAI
\cite{OpenAI}. It is used only in the final compact-superlevel
argument. The companion paper \cite{BanuelosSmooth2026} gives a
detailed exposition of its proof, together with related kernel results
and probabilistic interpretations. The principal arguments, techniques,
and proof structure in that exposition follow the OpenAI preprint.

\begin{smoothinput}[OpenAI {\cite{OpenAI}}]
\label{thm:smooth}
Let $G\subset\R^2$ be a bounded simply connected $C^\infty$ domain.
Then
\[
 \nabla u(x)\ne0\qquad(x\in G,\ 0\ne u\in E_1(G)).
\]
In particular, $G$ has property (HS1).
\end{smoothinput}

Our approximation uses only the case $\dim E_1(G)=1$ of
Theorem~\ref{thm:smooth}. The required input is the absence of interior
local extrema, which follows from the nonvanishing gradient.
Boundary attainment alone would not suffice for this argument.

The main contribution is an approximation of each prescribed
eigenfunction, including when $\dim E_1(\Omega)=2$. The attachment
point of a thin strip determines the limiting eigendirection. A
conformal exhaustion then gives smooth approximating domains.
The final passage from approximation to (HS1) uses the critical-set
argument of Judge and Mondal
\cite[Proposition~14, p.~160]{JudgeMondalCritical2025}.

Ba\~nuelos and Burdzy \cite{BB} introduced coupling arguments and
proved hot spots properties for obtuse triangles and for sufficiently
long convex planar domains with a line of symmetry. Jerison and
Nadirashvili \cite{JN} proved (HS1) for convex planar domains with two
orthogonal symmetry axes, and
Pascu \cite{Pascu2002} treated antisymmetric eigenfunctions on convex
$C^{1,\alpha}$ domains with one symmetry axis.

Atar and Burdzy
\cite{AB} proved (HS1) for \emph{lip domains}, that is, planar domains
bounded by the graphs of two Lipschitz functions with Lipschitz
constant at most one, and established simplicity
except for squares. Judge and Mondal \cite{JudgeMondal} (see also the
erratum cited there) proved (HS1) for every Euclidean triangle.
Hatcher \cite{Hatcher2025} gave an elementary proof for L-shaped domains
and established hot spots results for several families assembled from
such domains.

Neither simple connectedness nor the restriction to the plane can be
dropped. Burdzy and Werner \cite{BW1999},
Bass and Burdzy \cite{BassBurdzy2000}, and Burdzy \cite{Burdzy2005}
constructed planar counterexamples with holes. Atar and Burdzy
\cite{AtarBurdzy2026} constructed a bounded starlike Lipschitz
counterexample in $\R^3$ with simple first positive eigenvalue, and
de Dios Pont \cite[Theorem~1.2]{deDiosPont2024} constructed smooth
bounded convex counterexamples in sufficiently high dimensions.

Deng, Jiang, and Yang \cite{DengJiangYang2026} proved that, on a bounded
convex planar domain, every nonzero $u\in E_1(\Omega)$ satisfies
$|u(x)|<\max_{\partial\Omega}|u|$ for all $x\in\Omega$. Their
November~2, 2026 version states that its arguments do not rely on the
methods of \cite{OpenAI} and \cite{BurdzyNote}. Their proof uses a gradient
estimate of Li--Yau type and a near-boundary rigidity theorem. For each
$u$, this maximum-modulus conclusion gives at least one of the two
identities in \eqref{eq:main-HS2}. It does not by itself give both.

The notes \cite{Banuelos2026} and \cite{BurdzyNote} independently used
conformal exhaustion to pass from a smooth-domain input to a conclusion
of type (HS2) for less regular simply connected domains with a simple
first positive eigenvalue. The former also treats unbounded domains
of finite area. Here the nodal and multiplicity results for Lipschitz
domains, together with the strip perturbation, remove the simplicity
assumption. Lemma~\ref{lem:strict-boundary-attainment} supplies the
strict conclusion.

The idea of deforming the domain is taken from Jerison and Nadirashvili
\cite{JN}. Thin appendages are classical singular perturbations of the
Neumann problem \cite{Gadylshin2005,Jimbo,Jimbo1993}. The construction
below uses the location of the attachment to select an arbitrary
direction in the first eigenspace. The two-nodal-domain argument is
classical and also appears in \cite[Lemmas~2.3 and~2.4]{OpenAI}.
The nodal geometry is related to the mixed-boundary viewpoint of
Ba\~nuelos, Pang, and Pascu \cite{BPP}.

\begin{figure}[htbp]
\centering
\begin{tikzpicture}[scale=1.2,x=1cm,y=1cm,font=\small,line join=round,
  boundary/.style={draw=black!65,line width=.65pt,fill=black!3},
  inner/.style={draw=blue!65!black,line width=.7pt,rounded corners=2pt}]
\begin{scope}
 \node[font=\small\bfseries] at (1.4,3.65) {(a) Strip attached at $p_j$};
 \path[boundary] (.1,.35)--(.1,2)--(.75,2.4)--(1.37,2.4)
 --(1.37,3.05)--(1.63,3.05)--(1.63,2.4)--(1.9,2.4)
 --(2.65,1.85)--(2.25,1.5)--(2.72,.55)--(2,0)--(.7,.15)--cycle;
 \fill[orange!18] (1.37,2.4) rectangle (1.63,3.05);
 \draw[orange!70!black,line width=.8pt] (1.37,2.4)--(1.37,3.05)
 --(1.63,3.05)--(1.63,2.4);
 \draw[densely dashed,black!50] (1.37,2.4)--(1.63,2.4);
 \draw[blue!65!black,dashed,line width=.7pt] (1.05,2.4)
 .. controls (1.05,1.9) and (1.402,.525) .. (1.35,.075);
 \fill (1.05,2.4) circle (1.2pt);
 \fill (1.5,2.4) circle (1.2pt);
 \node[above left,inner sep=2pt] at (1.05,2.4) {$p_0$};
 \node[below right,inner sep=2pt] at (1.5,2.4) {$p_j$};
 \node[right,inner sep=1.5pt,font=\footnotesize] at (1.38,1.2) {$\{v=0\}$};
 \draw[black!55,line width=.3pt] (1.23,1.2)--(1.38,1.2);
 \node at (1.95,.5) {$\Omega$};
 \draw[<->] (1.37,3.24)--(1.63,3.24);
 \node[above,inner sep=1pt] at (1.5,3.24) {$\varepsilon_j$};
 \draw[<->] (1.9,2.45)--(1.9,3.05);
 \node[right,inner sep=2pt] at (1.9,2.77) {$h_j$};
 \node at (1.4,-.48) {length $h_j$, width $\varepsilon_j$};
\end{scope}
\begin{scope}[xshift=4.35cm]
 \node[font=\small\bfseries] at (1.4,3.65) {(b) Smooth exhaustion};
 \path[boundary] (.1,.35)--(.1,2)--(.75,2.4)--(1.37,2.4)
 --(1.37,3.05)--(1.63,3.05)--(1.63,2.4)--(1.9,2.4)
 --(2.65,1.85)--(2.25,1.5)--(2.72,.55)--(2,0)--(.7,.15)--cycle;
 \draw[inner,rounded corners=0pt] plot[smooth cycle,tension=.75] coordinates
 {(.32,.62) (.3,1.82) (.85,2.2) (1.5,2.3) (1.95,2.18) (2.38,1.84)
  (2.08,1.5) (2.42,.66) (1.93,.24) (.82,.33)};
 \draw[inner,dashed,rounded corners=0pt] plot[smooth cycle,tension=.75] coordinates
 {(.55,.8) (.52,1.68) (.98,1.97) (1.6,2.02) (2.0,1.8) (1.8,1.45)
  (2.1,.78) (1.8,.5) (.9,.55)};
 \node (Gjr) at (1.3,1.2) {$G_{j,r}$};
 \draw[black!55,line width=.3pt,->] (Gjr.east)--(2.2,1.2);
 \draw[->,blue!65!black] (.98,1.57)--(.98,2.22);
 \node[right,inner sep=2pt] at (.98,1.75) {$r\uparrow1$};
 \node[right,inner sep=2pt] at (1.65,2.76) {$\Omega_j$};
 \node at (1.4,-.48) {$G_{j,r}=F_j(r\D)$, fixed $j$};
\end{scope}
\begin{scope}[xshift=8.7cm]
 \node[font=\small\bfseries] at (1.4,3.65) {(c) Compact superlevel set};
 \draw[blue!65!black,line width=.8pt,fill=blue!2]
 (.28,.52) .. controls (.1,1) and (.2,1.85) .. (.87,2.16)
 .. controls (1.06,2.32) and (1.37,2.07) .. (1.44,2.42)
 .. controls (1.47,2.98) and (1.58,2.98) .. (1.61,2.42)
 .. controls (1.65,2.08) and (2.14,2.24) .. (2.41,1.85)
 .. controls (2.5,1.7) and (1.96,1.65) .. (2.1,1.4)
 .. controls (2.27,1.06) and (2.65,.63) .. (2.27,.38)
 .. controls (1.8,.02) and (.64,.1) .. (.28,.52)--cycle;
 \draw[fill=orange!10,draw=orange!65!black,line width=.8pt]
 (1.32,1.13) ellipse (.65 and .57);
 \node at (.95,1.35) {$U$};
 \fill (1.46,1.11) circle (1.4pt);
 \node[above right,inner sep=2pt] at (1.46,1.11) {$x_0$};
 \node[below,inner sep=2pt] at (1.32,.56) {$u=t$};
 \node at (1.25,1.98) {$G_{j,r}$};
 \node at (1.4,-.48) {$\overline U\subset G_{j,r}$};
\end{scope}
\end{tikzpicture}
\caption{The two domain approximations and the compact-superlevel argument,
shown schematically. In (a), a strip of width $\varepsilon_j$ and length
$h_j$ is attached at a point $p_j\in\Sigma_*$ near an endpoint $p_0$ of
the nodal arc of $v$. The dashed curve is the nodal arc of $v$. The
width is chosen after the length.
In (b), the interior curves represent conformal images of concentric
circles. The label $G_{j,r}$ refers to the domain bounded by the solid
curve.
In (c), $u=t$ on $\partial U$ while $u(x_0)>t$. A function that is
uniformly close to $u$ on $\overline U$ has a local maximum in $U$.}
\label{fig:domain-deformation}
\end{figure}
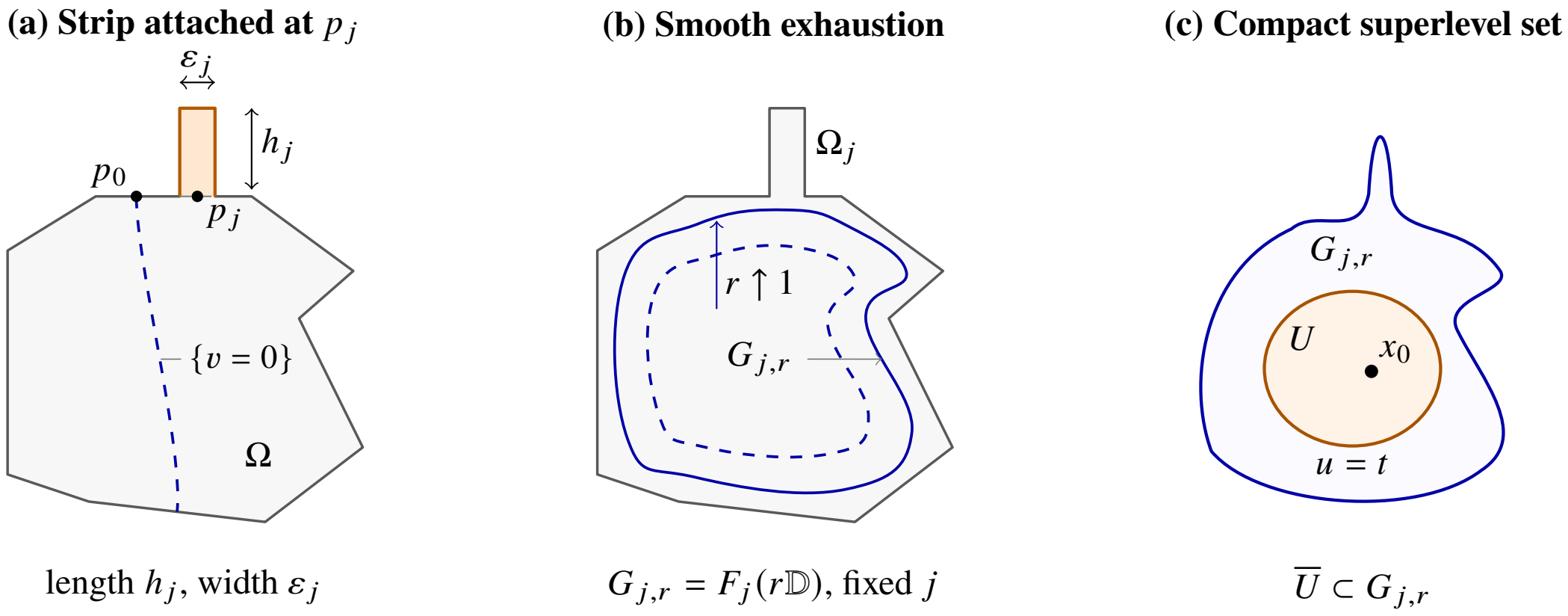

\paragraph{Outline of the proof of Theorem~\ref{thm:lipschitz}.}
The proof has two approximation steps, illustrated in
Figure~\ref{fig:domain-deformation}, followed by an argument with a
compact superlevel set and a critical-set argument that makes the
conclusion strict.

\emph{Step 1. Approximate $u$ by first eigenfunctions of domains whose
first positive eigenvalue is simple.}
Fix $u\in E_1(\Omega)$ with $\|u\|_{L^2(\Omega)}=1$.
Corollary~\ref{lem:lipschitz-multiplicity} gives
$\dim E_1(\Omega)\le2$. If $\dim E_1(\Omega)=1$, no approximation is
needed: we take $\Omega_j=\Omega$ and $u_j=u$. If
$\dim E_1(\Omega)=2$, choose $v\in E_1(\Omega)$ with
$\|v\|_{L^2(\Omega)}=1$ that is orthogonal to $u$ in $L^2(\Omega)$.
By Lemma~\ref{lem:nodal-arc}, the zero set of $v$ in
$\overline\Omega$ is an arc with two endpoints on $\partial\Omega$.
Corollary~\ref{lem:lipschitz-multiplicity} gives $u(p_0)\ne0$ at each
endpoint $p_0$. Attach to $\Omega$ a thin rectangular strip of width
$\varepsilon$ and length $h$ at a boundary point $p$. For $p$ in a dense
set $\Sigma_*\subset\partial\Omega$ whose complement has arclength zero
(Lemma~\ref{lem:strip-good-points}) and for $h$ small enough, the first
positive eigenvalue of the enlarged domain is simple for all small
$\varepsilon$, and as $\varepsilon\to0$ its normalized eigenfunction
converges in $H^1(\Omega)$, up to sign, to
\[
 \frac{u(p)u+v(p)v}{\sqrt{u(p)^2+v(p)^2}}
\]
(Proposition~\ref{prop:strip-evaluation}). At $p=p_0$ this function
equals $\pm u$, because $v(p_0)=0$. The point $p_0$ need not belong to
$\Sigma_*$. We therefore choose $p_j\in\Sigma_*$ with $p_j\to p_0$, then
lengths $h_j$, and then widths $\varepsilon_j$. This gives bounded
simply connected Lipschitz domains $\Omega_j\supset\Omega$ with simple
first positive Neumann eigenvalue and normalized eigenfunctions
$u_j\to u$ in $C^\infty_{\mathrm{loc}}(\Omega)$
(Proposition~\ref{prop:prescribed}).

\emph{Step 2. Exhaust each domain $\Omega_j$ by smooth domains.}
For fixed $j$, a conformal map $F_j$ of the unit disk $\D$ onto
$\Omega_j$ gives
\[
 G_{j,r}=F_j(r\D)\uparrow\Omega_j\qquad(r\uparrow1).
\]
These are $C^\infty$ domains. For $r$ close to $1$ their first
positive eigenvalue is simple, and their normalized first eigenfunctions
converge to $u_j$ locally uniformly in $\Omega_j$ as $r\uparrow1$
(Lemma~\ref{lem:eventual-simplicity}).

\emph{Conclusion.}
Suppose that, for some $t\in\R$, a connected component $U$ of
$\{u>t\}$ has compact closure in $\Omega$. Then $u=t$ on
$\partial U$. Fix $x_0\in U$. If $\tilde u$ is sufficiently close to
$u$ uniformly on $\overline U$, then
$\tilde u(x_0)>\max_{\partial U}\tilde u$, so $\tilde u$ has a local
maximum in $U$. By Steps~1 and~2 we may take for $\tilde u$ a first
eigenfunction of a smooth domain $G_{j,r}\supset\overline U$ with simple
first positive eigenvalue. This contradicts Theorem~\ref{thm:smooth}.
Applying the argument to $-u$ also excludes compactly contained
sublevel components. In particular, $u$ has no strict interior local
maximum or minimum. This argument uses the absence of interior local
extrema on the smooth domains. Boundary attainment of global extrema
alone would not suffice.

If $u$ attained its global maximum at
an interior point, real analyticity and the equation
$\Delta u=-\Lambda_\Omega^{\mathrm N}u$ would make that point either
an isolated maximum or a point of a real-analytic arc of maxima.
The first possibility is already excluded. The second is excluded by
the argument of Judge and Mondal cited above, as detailed in
Lemma~\ref{lem:strict-boundary-attainment}. Applying this to $-u$
gives (HS1).

\paragraph{Organization of the paper.}
Section~\ref{sec:lipschitz-prelim} proves the spectral, nodal, and
multiplicity facts used in the approximation.
Section~\ref{sec:lipschitz-bridge} constructs the strip perturbation,
computes the splitting of a double eigenvalue, and proves convergence
to a prescribed eigenfunction.
Section~\ref{sec:strict-attainment} combines conformal exhaustion with
Theorem~\ref{thm:smooth} and proves (HS1).
Section~\ref{sec:polygonal-comments} gives the simpler polygonal
reduction and explicit nonconvex starlike examples.

\FloatBarrier
\subsection{Notation and conventions}
\label{subsec:conventions}

\emph{Functions and domains.}
All eigenfunctions and test functions are real-valued, and all Sobolev
spaces are real. A \emph{domain} is a nonempty connected open subset of
$\R^2$, and $B(p,r)$ is the open disk of radius $r$ about $p$.
A bounded domain $\Omega\subset\R^2$ is called a \emph{Lipschitz domain}
if for every $p\in\partial\Omega$ there are a rigid motion that takes
$p$ to the origin, numbers $r,s>0$ and $L<\infty$, and a function
$\psi:(-r,r)\to(-s,s)$ with $\psi(0)=0$ and
$|\psi(x)-\psi(x')|\le L|x-x'|$, such that in the new coordinates
\[
 \Omega\cap Q
 =\{(x,y)\in Q:y<\psi(x)\},
 \qquad Q=(-r,r)\times(-s,s).
\]
We call $\Omega$ a \emph{$C^\infty$ domain} (respectively a
\emph{$C^{1,\alpha}$ domain}) if the functions $\psi$ in this
definition can be chosen in $C^\infty((-r,r))$ (respectively in
$C^{1,\alpha}((-r,r))$).

\emph{Eigenvalues.}
The Neumann eigenvalues of a bounded Lipschitz domain $\Omega$ are
indexed from zero and counted with multiplicity:
\[
 0=\Lambda_{\Omega,0}^{\mathrm N}<\Lambda_{\Omega,1}^{\mathrm N}
    \le\Lambda_{\Omega,2}^{\mathrm N}\le\cdots,
 \qquad
 \Lambda_\Omega^{\mathrm N}=\Lambda_{\Omega,1}^{\mathrm N}.
\]
We denote the lowest Dirichlet eigenvalue by
$\lambda_\Omega^{\mathrm D}$. The space $E_1(\Omega)$ is defined at the
beginning of Section~\ref{sec:introduction}. A \emph{first
eigenfunction} of $\Omega$ is a nonzero element of $E_1(\Omega)$. The
eigenvalue $\Lambda_\Omega^{\mathrm N}$ is \emph{simple} if
$\dim E_1(\Omega)=1$. The \emph{nodal set} of a function $u$ is its
zero set, and a \emph{nodal domain} of $u$ is a connected component of
$\{x\in\Omega:u(x)\ne0\}$. A \emph{normalized} eigenfunction has $L^2$
norm one on its domain.

\emph{Conformal maps.}
We write $\D=\{z\in\mathbb C:|z|<1\}$ and $\Sone=\partial\D$.
The letter $F$, with or without an index, denotes a conformal map
from $\D$ onto a simply connected domain. The symbol $ds$ denotes
arclength on the relevant boundary.

\section{Spectral preliminaries and nodal structure}
\label{sec:lipschitz-prelim}

The results of this section and the approximation results in
Section~\ref{sec:lipschitz-bridge} do not use
Theorem~\ref{thm:smooth}.

The notation and conventions are those of
Section~\ref{subsec:conventions}. The eigenvalues
$\Lambda_\Omega^{\mathrm N}$ and $\lambda_\Omega^{\mathrm D}$ have the
variational characterizations
\begin{align}
 \Lambda_\Omega^{\mathrm N}&=\inf_{\substack{v\in H^1(\Omega)\setminus\{0\}\\\int_\Omega v\,dx=0}}
       \frac{\int_\Omega|\nabla v|^2\,dx}{\int_\Omega|v|^2\,dx},
       \label{eq:neumann-minmax}\\
 \lambda_\Omega^{\mathrm D}&=\inf_{v\in H^1_0(\Omega)\setminus\{0\}}
       \frac{\int_\Omega|\nabla v|^2\,dx}{\int_\Omega|v|^2\,dx}.
       \label{eq:dirichlet-minmax}
\end{align}
We recall a special case of the theorem on p.~414 of Filonov's paper
\cite{Filonov}. If $\Omega\subset\R^n$, $n\ge2$,
has finite volume and $H^1(\Omega)\hookrightarrow L^2(\Omega)$ is
compact, then
\begin{equation}\label{eq:strict-spectral-gap}
 \Lambda_\Omega^{\mathrm N}<\lambda_\Omega^{\mathrm D}.
\end{equation}
For the history of this comparison, including the planar results of
P\'olya \cite{Polya1952} and Szeg\H{o} \cite{Szego1954}, see
\cite[p.~413]{Filonov}.
Neither simple connectedness nor eigenvalue simplicity is needed for
\eqref{eq:strict-spectral-gap}. We use it in the proof of
Lemma~\ref{lem:two-nodal-domains} to show that every nonzero
$u\in E_1(\Omega)$ takes positive and negative values on
$\partial\Omega$.

Let $\Omega$ be a bounded simply connected Lipschitz domain. Its
boundary is compact, connected (because $\Omega$ is simply connected),
and locally a Lipschitz graph; it is therefore a Jordan curve $J$. A
compactness argument with the local graphs shows that $J$ satisfies
Ahlfors' three-point condition
$\operatorname{diam}J(a,b)\le M|a-b|$ for $a,b\in J$, where $J(a,b)$ is
the arc of smaller diameter between $a$ and $b$. Hence $\Omega$ is a
$K$-quasidisk for some $K>1$, that is, the image of $\D$ under a
$K$-quasiconformal homeomorphism of the plane
\cite[Chapter~5]{Pommerenke}. By \cite[Proposition~5.3]{BGU2016},
applied with the unit disk as the first domain, every conformal map
$F:\D\to\Omega$ satisfies
\begin{equation}\label{eq:lipschitz-conformal-integrability}
 F'\in L^{2p}(\D),\qquad \rho=|F'|^2\in L^p(\D)
 \quad\text{for some }1<p<\infty,
\end{equation}
namely for every $1<p<K^2/(K^2-1)$. We use
\eqref{eq:lipschitz-conformal-integrability} in the proofs of
Lemmas~\ref{lem:conformal-reflection} and~\ref{lem:conformal-forms}.

\subsection{Nodal structure and multiplicity}
\label{subsec:lipschitz-multiplicity}

In what follows, we need several lemmas to determine the nodal set and
the boundary information used to select a prescribed eigendirection.
They will also give the multiplicity bound on every bounded simply
connected planar Lipschitz domain. The argument extends the
classical smooth-boundary conclusion of Nadirashvili
\cite{Nadirashvili1986}. See also \cite[Proposition~2.5]{BB}.
We transfer $u$ to the unit disk by a conformal map, reflect it to a
function on the sphere, and apply Kokarev's description
\cite{Kokarev2014} of the zero set of solutions of Schr\"odinger
equations with potential in $L^p$, $p>1$.

We begin with the two sign domains. The fact that both signs occur on
the boundary will allow each sign domain to remain connected after
reflection.

\begin{lemma}[Two nodal domains and boundary signs]
\label{lem:two-nodal-domains}
Let $\Omega\subset\R^2$ be a bounded simply connected Lipschitz domain,
and let $0\ne u\in E_1(\Omega)$. Then
$u\in C^\alpha(\overline\Omega)$ for some $\alpha>0$, the sets
$\{u>0\}$ and $\{u<0\}$ are its two nodal domains, and
\[
 \min_{\partial\Omega}u<0<\max_{\partial\Omega}u.
\]
\end{lemma}

\begin{proof}
Since $-\Delta u=\Lambda_\Omega^{\mathrm N}u\in L^2(\Omega)$
with homogeneous Neumann data, boundary regularity
\cite[Theorem~3.14(ii)]{Nittka2011} gives
$u\in C^\alpha(\overline\Omega)$ for some $\alpha>0$.

Since $\int_\Omega u=0$, there are at least two nodal domains. Suppose
there are at least three, and choose a nodal domain $\Omega_1$ on which
$u>0$ and a nodal domain $\Omega_2$ on which $u<0$. Put
$u_i=u\mathbf1_{\Omega_i}$. For $c>0$ the function
$(u-c)^+\mathbf1_{\Omega_1}$ vanishes in a neighborhood in $\Omega$ of
$\partial\Omega_1\cap\Omega$, because $u$ is continuous and $u=0$ there.
Hence it belongs to $H^1(\Omega)$, with gradient
$\mathbf1_{\Omega_1\cap\{u>c\}}\nabla u$. Letting $c\downarrow0$, and
arguing in the same way for $-u$ on $\Omega_2$, gives
$u_i\in H^1(\Omega)$ with $\nabla u_i=\mathbf1_{\Omega_i}\nabla u$.
Testing the weak eigenfunction equation with $u_i$ yields
\[
 \int_\Omega|\nabla u_i|^2=\Lambda_\Omega^{\mathrm N}\int_\Omega u_i^2,
 \qquad i=1,2.
\]
For
\[
 v=\left(\int_\Omega u_2\right)u_1
       -\left(\int_\Omega u_1\right)u_2\ne0
\]
the disjoint supports imply
\[
 \int_\Omega v=0,\qquad
 \frac{\int_\Omega|\nabla v|^2}{\int_\Omega v^2}=\Lambda_\Omega^{\mathrm N}.
\]
Thus $v\in E_1(\Omega)$ by the variational principle
\eqref{eq:neumann-minmax}, whereas $v=0$
on the third nodal domain, contradicting unique continuation.
The two nodal domains are therefore
\[
 \Omega_+=\{u>0\},\qquad \Omega_-=\{u<0\}.
\]
If $u\ge0$ on $\partial\Omega$, then $0\ne u^-=\max\{-u,0\}
\in H^1_0(\Omega)$ and the weak eigenfunction equation yields
\[
 \lambda_\Omega^{\mathrm D}
 \le\frac{\int_\Omega|\nabla u^-|^2}{\int_\Omega(u^-)^2}
 =\Lambda_\Omega^{\mathrm N},
\]
contrary to \eqref{eq:strict-spectral-gap}. Applying this to $-u$
also gives
\[
 \min_{\partial\Omega}u<0<\max_{\partial\Omega}u.
\]

\end{proof}

We now transfer the equation to the sphere, where the nodal set can be
studied without a boundary. For the remainder of this subsection, fix
a bounded simply connected
Lipschitz domain $\Omega$ and a conformal map $F:\D\to\Omega$.
Since $\Omega$ is a Jordan domain, Carath\'eodory's theorem extends
$F$ to a homeomorphism of the closures. By
\eqref{eq:lipschitz-conformal-integrability}, we may fix $p>1$ such that
\[
 \rho=|F'|^2\in L^p(\D).
\]
For $u\in E_1(\Omega)$, write $w=u\circ F$.

Identify $\mathbb S^2$ with $\mathbb C\cup\{\infty\}$ by stereographic
projection, so that $\D$ and $\{|z|>1\}\cup\{\infty\}$ are the two open
hemispheres and $\partial\D$ is the equator. The round metric and the
reflection in the equator are
\[
 g=\eta(z)|dz|^2,\qquad
 \eta(z)=\frac4{(1+|z|^2)^2},\qquad Jz=\frac1{\overline z}.
\]
Set
\[
 \widehat w(z)=
 \begin{cases}w(z),&|z|\le1,\\w(Jz),&|z|>1,\end{cases}
 \qquad
 \widehat\rho(z)=
 \begin{cases}\rho(z),&|z|<1,\\|z|^{-4}\rho(Jz),&|z|>1.\end{cases}
\]

\begin{lemma}[Conformal reflection]
\label{lem:conformal-reflection}
For every nonzero $u\in E_1(\Omega)$, the reflected function
$\widehat w$ belongs to $H^1(\mathbb S^2)\cap C(\mathbb S^2)$,
has exactly two nodal domains, and satisfies
\begin{equation}\label{eq:lipschitz-sphere-Schrodinger}
 (-\Delta_g+V)\widehat w=0,
 \qquad V=-\Lambda_\Omega^{\mathrm N}\widehat\rho/\eta\in L^p(\mathbb S^2),
\end{equation}
where the equation holds weakly on the whole sphere.
\end{lemma}

\begin{proof}
By Lemma~\ref{lem:two-nodal-domains}, $u$ is continuous on
$\overline\Omega$ and therefore bounded. 
Conformal invariance gives
\[
 w=u\circ F\in H^1(\D)\cap C(\overline\D),\qquad
 \int_\D|\nabla w|^2=\int_\Omega|\nabla u|^2,
\]
and
\begin{equation}\label{eq:lipschitz-weighted-Neumann}
 \int_\D\nabla w\cdot\nabla\varphi
 =\Lambda_\Omega^{\mathrm N}\int_\D\rho w\varphi,
 \qquad \varphi\in H^1(\D).
\end{equation}
Indeed, for $\varphi\in C^\infty(\overline\D)$ the function
$\varphi\circ F^{-1}$ belongs to $H^1(\Omega)$, and
\eqref{eq:lipschitz-weighted-Neumann} is the weak eigenfunction
equation for $u$ tested with it. The space $C^\infty(\overline\D)$ is
dense in $H^1(\D)$, and both sides of
\eqref{eq:lipschitz-weighted-Neumann} are continuous in $\varphi$ for
the $H^1(\D)$ norm, because, with $p'=p/(p-1)$,
\[
 \left|\int_\D\rho w\varphi\right|
 \le\|\rho\|_{L^p}\|w\|_\infty\|\varphi\|_{L^{p'}}
 \le C\|\varphi\|_{H^1(\D)}.
\]

Matching traces give
$\widehat w\in H^1(\mathbb S^2)\cap C(\mathbb S^2)$.
For $\Psi\in C^\infty(\mathbb S^2)$, use
$\varphi=(\Psi+\Psi\circ J)|_\D$ in
\eqref{eq:lipschitz-weighted-Neumann}. Change of variables on the
second hemisphere gives
\[
 \int_{\mathbb S^2}\langle\nabla_g\widehat w,\nabla_g\Psi\rangle_g\,dV_g
 =\Lambda_\Omega^{\mathrm N}\int_{\mathbb S^2}\frac{\widehat\rho}{\eta}
                  \widehat w\Psi\,dV_g.
\]
This proves \eqref{eq:lipschitz-sphere-Schrodinger} weakly across the equator.
Reflection gives the estimate
\[
 \int_{\mathbb S^2}|V|^p\,dV_g
 =2(\Lambda_\Omega^{\mathrm N})^p\int_\D\rho^p\eta^{1-p}
 \le C\int_\D\rho^p<\infty.
\]

By Lemma~\ref{lem:two-nodal-domains}, $\{w>0\}$ is a connected subset of $\D$, and its
image under $J$ is a connected subset of the opposite open hemisphere.
Every point $z\in\partial\D$ with $\widehat w(z)>0$ is the center of an
open disk on which $\widehat w>0$. This disk is connected, and it meets
both $\{w>0\}$ and its image under $J$. At least one such point exists,
because $\max_{\partial\Omega}u>0$. Hence $\{\widehat w>0\}$ is
connected. The same argument applies to $\{\widehat w<0\}$.

\end{proof}

The reflected equation allows us to apply Kokarev's nodal-structure
theorem. Since there are only two nodal domains, Euler's formula will
force the reflected nodal set to be a single regular curve.

\begin{lemma}[Nodal arc and boundary zeros]
\label{lem:nodal-arc}
For every nonzero $u\in E_1(\Omega)$, its nodal set in
$\overline\Omega$ is a single simple arc, with two distinct endpoints
on $\partial\Omega$ and all its other points in $\Omega$.

Let $z_0\in\partial\D$ satisfy $w(z_0)=0$. In a conformal coordinate
$\xi=t+is$ given by a M\"obius map that sends $\D$ to the upper
half-plane and $z_0$ to $0$, the reflected function has the expansion
\begin{equation}\label{eq:lipschitz-boundary-expansion}
 \widehat w(t,s)=a(w)t+O(|\xi|^{1+\delta'}),
 \qquad a(w)\ne0,
\end{equation}
for some $\delta'\in(0,1)$.
\end{lemma}

\begin{proof}
Use the reflected equation of Lemma~\ref{lem:conformal-reflection}.
For $0<\delta<\min\{1,2-2/p\}$, H\"older's inequality gives
\[
 \sup_{x\in\mathbb S^2}
 \int_{B_g(x,r)}\frac{|V(y)|}{d_g(x,y)^\delta}\,dV_g(y)
 \le C\|V\|_{L^p}r^{\,2-2/p-\delta}\longrightarrow0\qquad(r\downarrow0).
\]
Thus $V\in K^{2,\delta}$, the potential class of
\cite{Kokarev2014}. Let $\Gamma=\{\widehat w=0\}$. We use the following
facts about the nonzero continuous weak solution $\widehat w$ of
\eqref{eq:lipschitz-sphere-Schrodinger}, taken from
\cite[Theorem~1.3, Propositions~2.2 and~2.3, and Lemma~2.6]{Kokarev2014}.
The nodal-structure theorem applies to solutions with finitely many
nodal domains, as explained after its statement. No assertion about
the eigenvalue index of $\widehat w$ on the sphere is needed.
Every $q\in\Gamma$ has a finite vanishing order
$\operatorname{ord}_q\widehat w\ge1$. The set $\Gamma_2$ of points of
order at least two is finite. The set $\Gamma\setminus\Gamma_2$ has
finitely many connected components, and near each of its points
$\Gamma$ is a $C^1$ arc. Each component of $\Gamma\setminus\Gamma_2$ is
therefore a connected one-dimensional $C^1$ manifold, hence a closed
curve or an open arc. An open arc is a closed subset of
$\mathbb S^2\setminus\Gamma_2$, so the set of limit points of each of
its two ends is a connected subset of the finite set $\Gamma_2$, that
is, a single point. For $q\in\Gamma_2$, let $\deg(q)$ be the number of
ends of open arcs that converge to $q$; this is the degree defined in
\cite[Section~2.2]{Kokarev2014}. Then
\[
 \deg(q)\ge2\operatorname{ord}_q\widehat w\ge4
 \qquad(q\in\Gamma_2).
\]

Thus $\Gamma$ is a finite graph embedded in $\mathbb S^2$. Its vertices
are the points of $\Gamma_2$, its edges are the open arcs, and the
remaining components are closed curves. Insert a vertex of degree two
on each closed curve. Let $n_V$, $n_E$, and $n_C$ be the numbers of vertices,
edges, and connected components of this graph, and let $n_F$ be the
number of connected components of $\mathbb S^2\setminus\Gamma$. Euler's
formula for a finite graph embedded in the sphere reads
$n_V-n_E+n_F=1+n_C$. The components of $\mathbb S^2\setminus\Gamma$ are
the nodal domains of $\widehat w$, so $n_F=2$ by Lemma~\ref{lem:conformal-reflection}.
Hence
\[
 2=1+n_C+n_E-n_V,\qquad
 2n_E=\sum_q\deg(q)\ge2n_V,\qquad n_C\ge1.
\]
Therefore
\[
 n_C=1,\qquad n_E=n_V,\qquad \deg(q)=2\ \text{for every vertex }q.
\]
Since $\deg(q)\ge4$ on $\Gamma_2$, the set $\Gamma_2$ is empty. The
nodal set $\Gamma$ is a simple closed $C^1$ curve and
$\operatorname{ord}_q\widehat w=1$ for every $q\in\Gamma$.

Moreover $J\Gamma=\Gamma$ and
$\Gamma\cap\partial\D\ne\varnothing$, because $\widehat w$ takes both
signs on $\partial\D$. If $J|_\Gamma$ were the identity, then
$\Gamma\subset\partial\D$, hence $\Gamma=\partial\D$, a
contradiction. A continuous involution of a topological circle that
is not the identity and has a fixed point $q$ has exactly two fixed
points. Indeed, parametrize the circle by $[0,1]$ with both endpoints
sent to $q$. The involution corresponds to a homeomorphism $\iota$ of
$[0,1]$ that fixes $0$ and $1$ or exchanges them, with
$\iota\circ\iota=\mathrm{id}$. If $\iota$ were increasing, then
$\iota(x)>x$ would give $x=\iota(\iota(x))>\iota(x)$, and $\iota(x)<x$
would give $x<\iota(x)$ in the same way; so $\iota$ would be the
identity. Hence $\iota$ is decreasing, and it has exactly one fixed
point in $(0,1)$. Consequently
\[
 \#(\Gamma\cap\partial\D)=2.
\]
These points divide $\Gamma$ into two arcs exchanged by $J$, one
in each open hemisphere. Applying $F$ to the closed arc in $\overline\D$
proves the assertion of the lemma about the nodal set of $u$.

Now fix $z_0\in\partial\D$ with $w(z_0)=0$. Let $\xi=t+is$ be the coordinate given by
a M\"obius map that sends $\D$ to the upper half-plane and $z_0$ to
$0$. In this coordinate $\partial\D=\{s=0\}$ and
$J\xi=\overline\xi$. The expansion of
\cite[Proposition~2.2]{Kokarev2014} gives, whenever $w(z_0)=0$,
\[
 \widehat w(t,s)=a(w)t+b(w)s+O(|\xi|^{1+\delta'}),
 \qquad 0<\delta'<\delta.
\]
Here we use that $\operatorname{ord}_{z_0}\widehat w=1$ for
$w\ne0$. The expansion of \cite[Proposition~2.2]{Kokarev2014} holds in
some smooth chart at $z_0$. The passage from that chart to the
coordinate $\xi$ is smooth with invertible derivative at $z_0$; it
maps a nonzero linear leading term to a nonzero linear leading term and
changes the remainder by $O(|\xi|^2)$, which is $O(|\xi|^{1+\delta'})$
because $\delta'<1$. Since $\widehat w$ is even in $s$,
\[
 b(w)=0,\qquad a(w)\ne0\quad(w\ne0).
\]
This gives \eqref{eq:lipschitz-boundary-expansion}.
\end{proof}

The boundary expansion now controls the dimension of the space of
eigenfunctions vanishing at a given boundary point. This gives both
the multiplicity bound and the nonvanishing needed in the strip
construction.

\begin{corollary}[Multiplicity and boundary evaluation]
\label{lem:lipschitz-multiplicity}
Let $\Omega\subset\R^2$ be a bounded simply connected Lipschitz domain.
Then
\begin{equation}\label{eq:lipschitz-multiplicity}
 \dim E_1(\Omega)\le2.
\end{equation}
If equality holds, no basis of $E_1(\Omega)$ has a common zero on
$\partial\Omega$.
\end{corollary}

\begin{proof}
Fix $z_0\in\partial\D$ and use the conformal boundary coordinate
$\xi=t+is$ of Lemma~\ref{lem:nodal-arc}.
Let $\mathcal W=\{u\circ F:u\in E_1(\Omega)\}$ and
$\mathcal W_{z_0}=\ker(\mathrm{ev}_{z_0}:\mathcal W\to\R)$,
where $\mathrm{ev}_{z_0}(w)=w(z_0)$. The map
\[
 a:\mathcal W_{z_0}\longrightarrow\R,
 \qquad a(w)=\lim_{t\to0}\frac{\widehat w(t,0)}t
\]
is linear and, by \eqref{eq:lipschitz-boundary-expansion}, injective. Thus
\[
 \dim\mathcal W_{z_0}\le1,\qquad
 \dim E_1(\Omega)=\dim\mathcal W
 =\dim\mathcal W_{z_0}+\operatorname{rank}\mathrm{ev}_{z_0}\le2.
\]
If $\dim\mathcal W=2$, the same identity gives
$\operatorname{rank}\mathrm{ev}_{z_0}=1$ for every
$z_0\in\partial\D$. Hence no basis of $\mathcal W$ has a common
zero on $\partial\D$, and, because $F$ is a homeomorphism of the
closures, no basis of $E_1(\Omega)$ has a common zero on
$\partial\Omega$.
\end{proof}

\section{Approximation of a prescribed eigenfunction}
\label{sec:lipschitz-bridge}

Our aim is to approximate each prescribed eigenfunction by first
eigenfunctions on Lipschitz domains with a simple first positive
eigenvalue. We need several lemmas to control the spectrum under a
thin-strip perturbation and determine the selected eigendirection.
We then choose the attachment points using the nodal information from
Section~\ref{subsec:lipschitz-multiplicity}, obtaining
Proposition~\ref{prop:prescribed}.

\subsection{Thin-strip perturbations}
\label{subsec:strip-geometry}

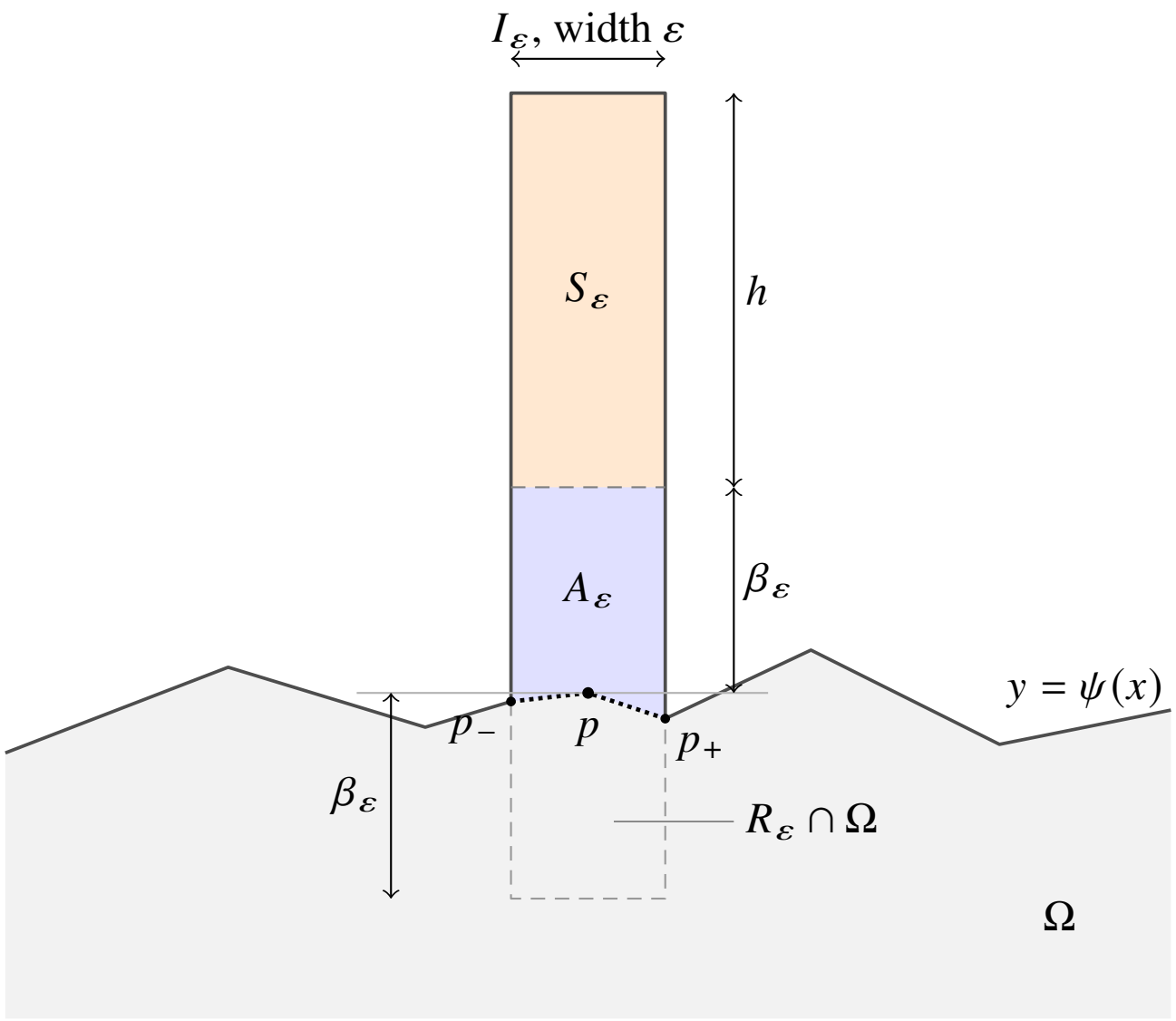
\begin{figure}[htbp]
\centering
\begin{tikzpicture}[x=1.25cm,y=1.25cm,font=\small,line join=round]
 \fill[black!5] (-3.4,-1.9)--(-3.4,-0.35)--(-2.1,0.15)--(-0.95,-0.2)
   --(-0.45,-0.05)--(0,0)--(0.45,-0.15)--(1.3,0.25)--(2.4,-0.3)
   --(3.4,-0.1)--(3.4,-1.9)--cycle;
 \fill[orange!18] (-0.45,1.2) rectangle (0.45,3.5);
 \fill[blue!12] (-0.45,-0.05)--(0,0)--(0.45,-0.15)--(0.45,1.2)
   --(-0.45,1.2)--cycle;
 \draw[black!70,line width=.7pt] (-3.4,-0.35)--(-2.1,0.15)--(-0.95,-0.2)
   --(-0.45,-0.05)--(-0.45,3.5)--(0.45,3.5)--(0.45,-0.15)--(1.3,0.25)
   --(2.4,-0.3)--(3.4,-0.1);
 \draw[black,densely dotted,line width=.9pt]
   (-0.45,-0.05)--(0,0)--(0.45,-0.15);
 \draw[black!50,densely dashed] (-0.45,1.2)--(0.45,1.2);
 \draw[black!40,densely dashed] (-0.45,-0.05)--(-0.45,-1.2)
   --(0.45,-1.2)--(0.45,-0.15);
 \draw[black!30] (-1.35,0)--(1.05,0);
 \fill (0,0) circle (1.2pt);
 \node[below,inner sep=5pt] at (0,0) {$p$};
 \fill (-0.45,-0.05) circle (1pt);
 \fill (0.45,-0.15) circle (1pt);
 \node[below left,inner sep=2pt] at (-0.45,-0.07) {$p_-$};
 \node[below right,inner sep=2pt] at (0.45,-0.17) {$p_+$};
 \node[right,inner sep=2pt] at (0.85,-0.75) {$R_\varepsilon\cap\Omega$};
 \draw[black!55,line width=.3pt] (0.85,-0.75)--(0.15,-0.75);
 \node at (0,0.6) {$A_\varepsilon$};
 \node at (0,2.35) {$S_\varepsilon$};
 \node at (2.75,-1.3) {$\Omega$};
 \node[above,inner sep=2pt] at (2.9,-0.15) {$y=\psi(x)$};
 \draw[<->] (-0.45,3.7)--(0.45,3.7);
 \node[above,inner sep=2pt] at (0,3.7) {$I_\varepsilon$, width $\varepsilon$};
 \draw[<->] (0.85,1.2)--(0.85,3.5);
 \node[right,inner sep=2pt] at (0.85,2.35) {$h$};
 \draw[<->] (0.85,0)--(0.85,1.2);
 \node[right,inner sep=2pt] at (0.85,0.65) {$\beta_\varepsilon$};
 \draw[<->] (-1.15,-1.2)--(-1.15,0);
 \node[left,inner sep=2pt] at (-1.15,-0.6) {$\beta_\varepsilon$};
\end{tikzpicture}
\caption{The strip attached at $p=0$. Up to the dotted graph arc and the
dashed segment $y=\beta_\varepsilon$, the rectangle
$R_\varepsilon=I_\varepsilon\times(-\beta_\varepsilon,\beta_\varepsilon+h)$
is the union of its part $R_\varepsilon\cap\Omega$ inside $\Omega$
(below the dotted graph arc),
the entrance region $A_\varepsilon$ between the graph and the level
$y=\beta_\varepsilon$, and the straight strip $S_\varepsilon$ of length
$h$. The thin horizontal line is the level $y=0$ of $p$, and $p_\pm$
are the endpoints of the graph arc. The solid curve is part of
$\partial\Omega_\varepsilon$.}
\label{fig:strip-geometry}
\end{figure}

Let $\Omega\subset\R^2$ be a bounded simply connected Lipschitz domain,
and fix $p\in\partial\Omega$. After a rigid motion, $p=0$ and, in a
rectangle $Q$ about $0$,
\[
 \Omega\cap Q=\{(x,y)\in Q:y<\psi(x)\},\qquad
 \psi(0)=0,\qquad \operatorname{Lip}(\psi)\le L.
\]
Choose $\beta_0>L/2+1$ and fix $h>0$ sufficiently small that the
construction below lies in $Q$ for all sufficiently small
$\varepsilon>0$. Set
\begin{equation}\label{eq:strip-domain}
 \begin{split}
 I_\varepsilon&=(-\varepsilon/2,\varepsilon/2),
 \qquad \beta_\varepsilon=\beta_0\varepsilon,\\
 R_\varepsilon&=I_\varepsilon\times(-\beta_\varepsilon,\beta_\varepsilon+h),
 \qquad \Omega_\varepsilon=\Omega\cup R_\varepsilon.
 \end{split}
\end{equation}
The added part consists, up to sets of area zero, of the entrance
region and the straight strip
\[
 A_\varepsilon=\{(x,y):x\in I_\varepsilon,
                    \ \psi(x)<y<\beta_\varepsilon\},
 \qquad
 S_\varepsilon=I_\varepsilon\times(\beta_\varepsilon,\beta_\varepsilon+h).
\]
In particular,
\[
 |A_\varepsilon|=O(\varepsilon^2),\qquad
 |S_\varepsilon|=\varepsilon h,\qquad
 \Omega\subset\Omega_\varepsilon.
\]
The sets $I_\varepsilon$, $A_\varepsilon$, $S_\varepsilon$, and
$R_\varepsilon$ are shown in Figure~\ref{fig:strip-geometry}.

On $I_\varepsilon$ we have $|\psi|\le L\varepsilon/2<\beta_\varepsilon$.
Hence the open rectangle $R_\varepsilon$ meets $\Omega$ in the connected
set $\{(x,y):x\in I_\varepsilon,\ -\beta_\varepsilon<y<\psi(x)\}$, and
$\Omega_\varepsilon$ is the union of two simply connected domains with
connected intersection. It is therefore simply connected. Its boundary
is obtained from $\partial\Omega$ by replacing the graph arc over
$I_\varepsilon$ by the segments
$\{\pm\varepsilon/2\}\times[\psi(\pm\varepsilon/2),\beta_\varepsilon+h]$
of the vertical sides of $R_\varepsilon$ and by the horizontal top of
$R_\varepsilon$. These three segments meet $\partial\Omega$ only at the
two points $p_\pm=(\pm\varepsilon/2,\psi(\pm\varepsilon/2))$, so the new
boundary is a Jordan curve. Set $e_\pm=(L+1,\pm1)$ and
$n_\pm=(\mp1,L+1)$, so that $e_\pm\cdot n_\pm=0$. Traversed from left to
right, the boundary near $p_-$ consists of the graph followed by the
left side traversed upward, and near $p_+$ it consists of the right
side traversed downward followed by the graph. For $x_1<x_2$,
\[
 \bigl(x_2-x_1,\psi(x_2)-\psi(x_1)\bigr)\cdot e_\pm\ge x_2-x_1,\qquad
 (0,\pm1)\cdot e_\pm=1.
\]
Let $z,z'$ be two points of the boundary near $p_-$, with $z$ before
$z'$. If both lie on the graph or both lie on the left side, the
display gives $(z'-z)\cdot e_+\ge|z'-z|/\sqrt{1+L^2}$. If $z$ lies on
the graph and $z'$ on the side, the same inequality follows by adding
the inequalities for the pairs $(z,p_-)$ and $(p_-,z')$ and using
$|z'-z|\le|z'-p_-|+|p_--z|$. Hence, in the orthonormal coordinates
$(\xi,\eta)$ with origin $p_-$ and axes in the directions $e_+$ and
$n_+$, the boundary near $p_-$ is the graph $\eta=g(\xi)$ of a
Lipschitz function. Moreover, $\Omega_\varepsilon$ coincides near $p_-$
with $\{\eta<g(\xi)\}$: for a boundary point $z$ near $p_-$ and small
$s>0$, the point $z+sn_+$ lies to the left of the line
$x=-\varepsilon/2$ and above the graph of $\psi$, hence outside
$\Omega_\varepsilon$, and the point $z-sn_+$ lies in $R_\varepsilon$ or
below the graph of $\psi$, hence in $\Omega_\varepsilon$. The same
argument applies at $p_+$ with $e_-$ and $n_-$.
The top corners are right angles. Hence $\Omega_\varepsilon$ is bounded,
simply connected, and Lipschitz. The numbers $r$ and $s$ in the
definition of a Lipschitz domain, applied to $\Omega_\varepsilon$,
depend on $\varepsilon$.

In the estimates that follow, constants may depend on the fixed
domain, chart, $\beta_0$, and $h$, but not on $\varepsilon$. Write
\[
 a_\varepsilon(f,g)=\int_{\Omega_\varepsilon}\nabla f\cdot\nabla g,
 \qquad
 b_\varepsilon(f,g)=\int_{\Omega_\varepsilon}fg.
\]

We must first rule out concentration of low-energy eigenfunctions in
the added region. The next two lemmas control the entrance and reduce
the behavior in the strip to a one-dimensional problem.

\begin{lemma}[Entrance trace and mass estimates]\label{lem:strip-traces}
In the setting of \eqref{eq:strip-domain} there is a constant $C$,
independent of $\varepsilon$ and $\phi$, such that the following holds
for all sufficiently small $\varepsilon$ and every
$\phi\in H^1(\Omega_\varepsilon)$. We have
\begin{equation}\label{eq:strip-base-trace}
 \left|\int_{I_\varepsilon}\phi(x,\beta_\varepsilon)\,dx\right|
 \le C\varepsilon^{3/4}\|\phi\|_{H^1(\Omega)}
       +C\varepsilon\|\nabla\phi\|_{L^2(A_\varepsilon)},
\end{equation}
and
\begin{equation}\label{eq:strip-entrance-mass}
 \|\phi\|_{L^2(A_\varepsilon)}^2
 \le C\varepsilon^{3/2}\|\phi\|_{H^1(\Omega)}^2
       +C\varepsilon^2\|\nabla\phi\|_{L^2(A_\varepsilon)}^2.
\end{equation}
\end{lemma}

\begin{proof}
Put $\theta(x)=\operatorname{Tr}\phi(x,\psi(x))$ and
$\ell_\varepsilon(x)=\beta_\varepsilon-\psi(x)$. Since
$|\psi|\le L\varepsilon/2$ on $I_\varepsilon$ and
$\beta_\varepsilon=\beta_0\varepsilon$ with $\beta_0>L/2+1$, we have
$0<\ell_\varepsilon\le C\varepsilon$. The trace theorem
$H^1(\Omega)\to H^{1/2}(\partial\Omega)$ for the Lipschitz domain
$\Omega$ and the Sobolev embedding
$H^{1/2}(\partial\Omega)\subset L^4(\partial\Omega)$ give
\[
 \|\theta\|_{L^4(I_\varepsilon)}
 \le\|\operatorname{Tr}\phi\|_{L^4(\partial\Omega)}
 \le C\|\phi\|_{H^1(\Omega)}.
\]
The traces from $\Omega$ and $A_\varepsilon$ agree. Hence, for almost
every $x\in I_\varepsilon$,
\[
 \phi(x,y)=\theta(x)+\int_{\psi(x)}^y\partial_y\phi(x,s)\,ds,
 \qquad \psi(x)<y\le \beta_\varepsilon.
\]
H\"older's inequality therefore yields
\[
 \begin{aligned}
 \left|\int_{I_\varepsilon}\phi(x,\beta_\varepsilon)\,dx\right|
 &\le\varepsilon^{3/4}\|\theta\|_{L^4(I_\varepsilon)}
       +|A_\varepsilon|^{1/2}\|\partial_y\phi\|_{L^2(A_\varepsilon)}\\
 &\le C\varepsilon^{3/4}\|\phi\|_{H^1(\Omega)}
       +C\varepsilon\|\nabla\phi\|_{L^2(A_\varepsilon)},
 \end{aligned}
\]
which is \eqref{eq:strip-base-trace}. Similarly,
\[
 \begin{aligned}
 \int_{A_\varepsilon}|\phi|^2
 &\le 2\int_{I_\varepsilon}\ell_\varepsilon|\theta|^2
       +C\varepsilon^2\int_{A_\varepsilon}|\partial_y\phi|^2\\
 &\le C\varepsilon^{3/2}\|\theta\|_{L^4(I_\varepsilon)}^2
       +C\varepsilon^2\|\nabla\phi\|_{L^2(A_\varepsilon)}^2,
 \end{aligned}
\]
proving \eqref{eq:strip-entrance-mass}.
\end{proof}

We next separate the longitudinal average from the variation across
the strip. The entrance estimate will force this normalized average
to vanish at the base in the limit for families with bounded $H^1$ norms.

\begin{lemma}[Averaging in the strip]
\label{lem:strip-average}
In the setting of \eqref{eq:strip-domain}, let
$\phi\in H^1(\Omega_\varepsilon)$ and define
\begin{equation}\label{eq:strip-average}
 g_\varepsilon(t)=\varepsilon^{-1/2}
       \int_{I_\varepsilon}\phi(x,\beta_\varepsilon+t)\,dx,
       \qquad 0<t<h.
\end{equation}
Then $g_\varepsilon\in H^1(0,h)$,
\begin{equation}\label{eq:strip-average-bound}
 \begin{aligned}
 \|g_\varepsilon\|_{H^1(0,h)}
 &\le\|\phi\|_{H^1(S_\varepsilon)},\\
 |g_\varepsilon(0)|
 &\le C\varepsilon^{1/4}\|\phi\|_{H^1(\Omega)}
       +C\varepsilon^{1/2}\|\nabla\phi\|_{L^2(A_\varepsilon)},
 \end{aligned}
\end{equation}
and
\begin{equation}\label{eq:strip-transverse}
 \int_{S_\varepsilon}
 \left|\phi(x,y)-\varepsilon^{-1/2}
                    g_\varepsilon(y-\beta_\varepsilon)\right|^2\,dx\,dy
 \le C\varepsilon^2\int_{S_\varepsilon}|\partial_x\phi|^2.
\end{equation}
The constants are independent of $\varepsilon$ and $\phi$.
\end{lemma}

\begin{proof}
For $g_\varepsilon$ in \eqref{eq:strip-average}, Jensen's inequality
and \eqref{eq:strip-base-trace} give
\[
 \begin{aligned}
 \int_0^h(|g_\varepsilon|^2+|g_\varepsilon'|^2)
 &\le\int_{S_\varepsilon}(|\phi|^2+|\partial_y\phi|^2),\\
 |g_\varepsilon(0)|
 &\le C\varepsilon^{1/4}\|\phi\|_{H^1(\Omega)}
       +C\varepsilon^{1/2}\|\nabla\phi\|_{L^2(A_\varepsilon)}.
 \end{aligned}
\]
Finally, the mean of $\phi(\cdot,\beta_\varepsilon+t)$ is
$\varepsilon^{-1/2}g_\varepsilon(t)$. The interval Poincar\'e
inequality gives
\[
 \int_0^h\int_{I_\varepsilon}
 \big|\phi(x,\beta_\varepsilon+t)-\varepsilon^{-1/2}g_\varepsilon(t)\big|^2
 \,dx\,dt
 \le\frac{\varepsilon^2}{\pi^2}
       \int_{S_\varepsilon}|\partial_x\phi|^2,
\]
and hence \eqref{eq:strip-transverse}.
\end{proof}

For the next two lemmas, fix $T>0$ and assume that the strip length
$h$ in \eqref{eq:strip-domain} satisfies
\begin{equation}\label{eq:strip-threshold}
 T<\nu_h:=\frac{\pi^2}{4h^2}.
\end{equation}
The number $\nu_h$ is the first eigenvalue on $(0,h)$ with a Dirichlet
condition at the entrance and a Neumann condition at the tip. The
inequality \eqref{eq:strip-threshold} prevents concentration of
low-energy eigenfunctions in the strip.

\begin{lemma}[Compactness below the strip threshold]
\label{lem:strip-compactness}
Let $\Omega$, $p$, and $\Omega_\varepsilon$ be as in
\eqref{eq:strip-domain}, and assume \eqref{eq:strip-threshold}.
Suppose that $\phi_\varepsilon$ is a normalized Neumann eigenfunction
on $\Omega_\varepsilon$ with eigenvalue
$0\le\gamma_\varepsilon\le T$. Then
\[
 \|\phi_\varepsilon\|_{L^2(A_\varepsilon\cup S_\varepsilon)}
 \longrightarrow0.
\]
Along a subsequence as $\varepsilon\downarrow0$, there are
$\gamma\in[0,T]$ and $\phi\in H^1(\Omega)$ such that
\[
 \begin{gathered}
 \gamma_\varepsilon\longrightarrow\gamma,\qquad
 \|\phi\|_{L^2(\Omega)}=1,\\
 \phi_\varepsilon|_\Omega\rightharpoonup\phi
       \quad\text{in }H^1(\Omega),\qquad
 \phi_\varepsilon|_\Omega\longrightarrow\phi
       \quad\text{in }L^2(\Omega),\\
 \int_\Omega\nabla\phi\cdot\nabla f
       =\gamma\int_\Omega\phi f
       \quad(f\in H^1(\Omega)).
 \end{gathered}
\]
For finitely many such eigenfunction families $\phi_\varepsilon^1,\ldots,
\phi_\varepsilon^m$ that are orthonormal for each $\varepsilon$, the
subsequence can be chosen simultaneously, and the limiting
eigenfunctions are orthonormal in $L^2(\Omega)$.
\end{lemma}

\begin{proof}
Let $\phi_\varepsilon$ be a normalized Neumann eigenfunction on
$\Omega_\varepsilon$, with eigenvalue $\gamma_\varepsilon\le T$. Then
\[
 \|\phi_\varepsilon\|_{H^1(\Omega_\varepsilon)}^2
 =1+\gamma_\varepsilon\le1+T.
\]
Rellich compactness on $\Omega$ gives, along a subsequence,
\[
 \gamma_\varepsilon\to\gamma\in[0,T],\qquad
 \phi_\varepsilon|_\Omega\rightharpoonup\phi\text{ in }H^1(\Omega),
 \qquad
 \phi_\varepsilon|_\Omega\to\phi\text{ in }L^2(\Omega).
\]
Let $g_\varepsilon$ be the average \eqref{eq:strip-average} of
$\phi_\varepsilon$. Testing with
$\varepsilon^{-1/2}v(y-\beta_\varepsilon)$ on $S_\varepsilon$, extended
by zero, gives
\[
 \int_0^h g_\varepsilon'v'
 =\gamma_\varepsilon\int_0^h g_\varepsilon v
 \qquad\bigl(v\in H^1(0,h),\ v(0)=0\bigr).
\]
Thus
\[
 -g_\varepsilon''=\gamma_\varepsilon g_\varepsilon,
 \qquad g_\varepsilon'(h)=0,
\]
and, also when $\gamma_\varepsilon=0$,
\begin{equation}\label{eq:strip-average-profile}
 g_\varepsilon(t)=g_\varepsilon(0)
 \frac{\cos(\sqrt{\gamma_\varepsilon}(h-t))}
      {\cos(\sqrt{\gamma_\varepsilon}h)}.
\end{equation}
Since $\sqrt T h<\pi/2$, \eqref{eq:strip-average-bound} implies
\[
 \cos(\sqrt{\gamma_\varepsilon}h)\ge\cos(\sqrt T h)>0,
 \qquad
 \|g_\varepsilon\|_{H^1(0,h)}
 \le C_T|g_\varepsilon(0)|\longrightarrow0.
\]
Orthogonality of the sectional mean and mean-zero part, together
with \eqref{eq:strip-transverse}, yields
\[
 \begin{aligned}
 \|\phi_\varepsilon\|_{L^2(S_\varepsilon)}^2
 &=\|g_\varepsilon\|_{L^2(0,h)}^2
   +\int_{S_\varepsilon}
    |\phi_\varepsilon-\varepsilon^{-1/2}
                    g_\varepsilon(y-\beta_\varepsilon)|^2\\
 &\le\|g_\varepsilon\|_{L^2(0,h)}^2+C\varepsilon^2(1+T)
 \longrightarrow0.
 \end{aligned}
\]
By \eqref{eq:strip-entrance-mass},
\[
 \|\phi_\varepsilon\|_{L^2(A_\varepsilon)}^2
 \le C(\varepsilon^{3/2}+\varepsilon^2)(1+T)\longrightarrow0.
\]
Since $\|\phi_\varepsilon\|_{L^2(\Omega_\varepsilon)}=1$ and the $L^2$
mass of $\phi_\varepsilon$ on $A_\varepsilon\cup S_\varepsilon$ tends to
zero,
$\|\phi\|_{L^2(\Omega)}=\lim\|\phi_\varepsilon\|_{L^2(\Omega)}=1$.

For $f\in H^1(\Omega)$ vanishing near $p$, its extension by zero to
$\Omega_\varepsilon\setminus\Omega$ belongs to $H^1(\Omega_\varepsilon)$
for small $\varepsilon$. Passing to the limit gives
\[
 \int_\Omega\nabla\phi\cdot\nabla f=\gamma\int_\Omega\phi f.
\]
Such functions $f$ are dense in $H^1(\Omega)$. Indeed, for $0<r<1$, put
\[
 \eta_r(x)=
 \begin{cases}
 0,&|x-p|\le r^2,\\
 \dfrac{\log(|x-p|/r^2)}{\log(1/r)},&r^2<|x-p|<r,\\
 1,&|x-p|\ge r.
 \end{cases}
\]
Then
\[
 \int_{\R^2}|\nabla\eta_r|^2=\frac{2\pi}{\log(1/r)}\to0,
 \qquad \eta_rf\to f\text{ in }H^1(\Omega)
\]
for every $f$ smooth near $\overline\Omega$. Such restrictions are
dense by the Lipschitz extension theorem. Consequently the weak
eigenfunction identity holds for all $f\in H^1(\Omega)$.

The same argument for finitely many orthonormal eigenfunctions gives
\[
 \langle\phi^i,\phi^j\rangle_{L^2(\Omega)}
 =\lim_{\varepsilon\downarrow0}
   \langle\phi_\varepsilon^i,\phi_\varepsilon^j
                      \rangle_{L^2(\Omega_\varepsilon)}
 =\delta_{ij}.
\]
\end{proof}

The compactness lemma preserves the orthogonality of finitely many
eigenfunctions. Combined with the min--max principle, this controls
both the location and the number of eigenvalues near
$\Lambda_\Omega^{\mathrm N}$.

\begin{lemma}[Preservation of the two lowest positive eigenvalues]
\label{lem:strip-spectrum}
Let $\Omega$, $p$, and $\Omega_\varepsilon$ be as in
\eqref{eq:strip-domain}, and assume \eqref{eq:strip-threshold}.
Suppose that
\[
 0=\Lambda_{\Omega,0}^{\mathrm N}<\Lambda_{\Omega,1}^{\mathrm N}=\Lambda_{\Omega,2}^{\mathrm N}
 <T<\Lambda_{\Omega,3}^{\mathrm N}.
\]
Then
\begin{equation}\label{eq:strip-cluster-convergence}
 \begin{aligned}
 \Lambda_{\Omega_\varepsilon,1}^{\mathrm N},\Lambda_{\Omega_\varepsilon,2}^{\mathrm N}
       &\longrightarrow\Lambda_\Omega^{\mathrm N}\qquad(\varepsilon\downarrow0),\\
 \Lambda_{\Omega_\varepsilon,3}^{\mathrm N}&>T
       \quad\text{for all sufficiently small }\varepsilon.
 \end{aligned}
\end{equation}
In particular, for small $\varepsilon$ the distance from
$\{\Lambda_{\Omega_\varepsilon,1}^{\mathrm N},\Lambda_{\Omega_\varepsilon,2}^{\mathrm N}\}$
to the rest of the Neumann spectrum of $\Omega_\varepsilon$ is bounded
below by a positive constant independent of $\varepsilon$.
\end{lemma}

\begin{proof}
Given $\kappa>0$, approximate a basis of
$\operatorname{span}\{1\}\oplus E_1(\Omega)$ in $H^1(\Omega)$ by
$F_0,F_1,F_2\in C_c^\infty(\R^2)$ so that
\[
 \dim\operatorname{span}\{F_i|_\Omega:0\le i\le2\}=3,
 \qquad
 \max_{a\ne0}
 \frac{\int_\Omega|\nabla\sum_{i=0}^2a_iF_i|^2}
      {\int_\Omega|\sum_{i=0}^2a_iF_i|^2}
 \le\Lambda_\Omega^{\mathrm N}+\kappa.
\]
Since $|\Omega_\varepsilon\setminus\Omega|=O(\varepsilon)$,
\[
 \left(\int_{\Omega_\varepsilon}F_iF_j\right)_{i,j=0}^2
 \longrightarrow\left(\int_\Omega F_iF_j\right)_{i,j=0}^2>0,
 \qquad
 \int_{\Omega_\varepsilon}\nabla F_i\cdot\nabla F_j
 \longrightarrow\int_\Omega\nabla F_i\cdot\nabla F_j.
\]
Min--max therefore gives
$\limsup_{\varepsilon\downarrow0}\Lambda_{\Omega_\varepsilon,2}^{\mathrm N}
\le\Lambda_\Omega^{\mathrm N}+\kappa$ for every $\kappa>0$, and hence
\[
 \limsup_{\varepsilon\downarrow0}\Lambda_{\Omega_\varepsilon,1}^{\mathrm N}
 \le
 \limsup_{\varepsilon\downarrow0}\Lambda_{\Omega_\varepsilon,2}^{\mathrm N}
 \le\Lambda_\Omega^{\mathrm N}.
\]
Hence both eigenvalues are eventually below $T$. For either positive
eigenvalue, choose corresponding normalized eigenfunctions
$\phi_\varepsilon$. By Lemma~\ref{lem:strip-compactness}, every
subsequential limit is an eigenvalue of $\Omega$ with a normalized
limiting eigenfunction $\phi$. Moreover,
\[
 \int_\Omega\phi
 =\lim_{\varepsilon\downarrow0}\int_\Omega\phi_\varepsilon
 =-\lim_{\varepsilon\downarrow0}
       \int_{\Omega_\varepsilon\setminus\Omega}\phi_\varepsilon=0,
\]
where the last limit follows by Cauchy--Schwarz. Thus each limit is
at least $\Lambda_\Omega^{\mathrm N}$, and
\[
 \Lambda_\Omega^{\mathrm N}\le\liminf\Lambda_{\Omega_\varepsilon,i}^{\mathrm N}
 \le\limsup\Lambda_{\Omega_\varepsilon,i}^{\mathrm N}\le\Lambda_\Omega^{\mathrm N},
 \qquad i=1,2.
\]

If $\Lambda_{\Omega_\varepsilon,3}^{\mathrm N}\le T$ along a sequence,
apply Lemma~\ref{lem:strip-compactness} to the first four orthonormal
eigenfunctions. It gives orthonormal limits $\phi^0,\ldots,\phi^3$ on
$\Omega$, with eigenvalues at most $T$.
Their span $V$ therefore satisfies
\[
 \dim V=4,
 \qquad
 \Lambda_{\Omega,3}^{\mathrm N}\le
 \max_{0\ne f\in V}\frac{\int_\Omega|\nabla f|^2}{\int_\Omega f^2}
 \le T,
\]
contrary to $T<\Lambda_{\Omega,3}^{\mathrm N}$. This proves
\eqref{eq:strip-cluster-convergence}. The uniform gaps follow from
$\Lambda_\Omega^{\mathrm N}>0$ and $T>\Lambda_\Omega^{\mathrm N}$.
\end{proof}

\subsection{First-order splitting of a double eigenvalue}
\label{subsec:strip-splitting}

Assume $\Lambda^{\mathrm N}_{\Omega,1}=\Lambda^{\mathrm N}_{\Omega,2}$.
We compute $\Lambda^{\mathrm N}_{\Omega_\varepsilon,1}$ and
$\Lambda^{\mathrm N}_{\Omega_\varepsilon,2}$ up to $o(\varepsilon)$
(Proposition~\ref{prop:strip-evaluation}). The result involves the
values, at the point $p$ at which the strip is attached, of the
functions in $E_1(\Omega)$. We prove it for points $p$ of the following
kind.

\begin{definition}[Admissible points]
\label{def:admissible}
Let $\Omega$ be a bounded simply connected Lipschitz domain. We denote
by $\Sigma_*=\Sigma_*(\Omega)$ the set of points $p\in\partial\Omega$
such that, for every $u\in E_1(\Omega)$,
\begin{equation}\label{eq:strip-good-energy}
 \int_{\Omega\cap B(p,r)}|\nabla u|^2=o(r)
       \qquad(r\downarrow0).
\end{equation}
The points of $\Sigma_*$ are called \emph{admissible}.
\end{definition}

The attachment point will eventually approach a prescribed nodal
endpoint. We therefore need admissible points arbitrarily close to
every boundary point.

\begin{lemma}[Admissible boundary points]
\label{lem:strip-good-points}
Let $\Omega$ be a bounded simply connected Lipschitz domain with
$\dim E_1(\Omega)=2$. The set $\Sigma_*$ is dense in $\partial\Omega$,
and its complement in $\partial\Omega$ has arclength zero.
\end{lemma}

\begin{proof}
Choose an $L^2(\Omega)$-orthonormal basis $u_1,u_2$ of $E_1(\Omega)$.
Set
\[
 \mathcal M=\mathbf1_\Omega(|\nabla u_1|^2+|\nabla u_2|^2)\,dx,
 \qquad \sigma=\mathcal H^1|_{\partial\Omega}.
\]
The measure $\mathcal M$ is absolutely continuous with respect to
Lebesgue measure, and $\sigma$ is concentrated on $\partial\Omega$,
which has Lebesgue measure zero. Hence $\mathcal M\perp\sigma$. By the
differentiation theory for Radon measures
\cite[Section~1.6]{EvansGariepy} (the singular part of a Radon measure
has derivative zero almost everywhere with respect to the reference
measure), the derivative of $\mathcal M$ with
respect to $\sigma$ vanishes $\sigma$-almost everywhere. Since
$\partial\Omega$ is locally a Lipschitz graph, there are $c,C>0$ with
$cr\le\sigma(B(p,r))\le Cr$ for all $p\in\partial\Omega$ and all small
$r$. Hence, for $\sigma$-almost every $p$,
\[
 \frac{\mathcal M(B(p,r))}{\sigma(B(p,r))}\longrightarrow0,
 \qquad cr\le\sigma(B(p,r))\le Cr,
 \qquad \mathcal M(B(p,r))=o(r).
\]
Since every $u\in E_1(\Omega)$ is a linear combination of $u_1$ and
$u_2$, such a point $p$ satisfies \eqref{eq:strip-good-energy} for every
$u\in E_1(\Omega)$. Thus the complement of $\Sigma_*$ has arclength
zero. The set $\Sigma_*$ is dense, since every relatively open boundary
arc has positive $\sigma$ measure.
\end{proof}

Admissibility supplies the error scale needed for the first-order
calculation. We use it to make each eigenfunction constant near the
entrance, so that it can be joined to a one-dimensional strip profile.

\begin{lemma}[Local flattening]
\label{lem:strip-flattening}
Let $\Omega$ be as in Lemma~\ref{lem:strip-good-points}, let
$u_1,u_2$ be an $L^2(\Omega)$-orthonormal basis of $E_1(\Omega)$, and
fix $p\in\Sigma_*$. Consider the strip construction
\eqref{eq:strip-domain} at $p$. For all sufficiently small
$\varepsilon$ there are constants $c_i(\varepsilon)$ and functions
$\zeta_{i,\varepsilon}\in H^1(\Omega)$, $i=1,2$, such that
$u_i+\zeta_{i,\varepsilon}=c_i(\varepsilon)$ in
$\Omega\cap W_\varepsilon$ for some open set
$W_\varepsilon\subset\R^2$ containing the graph arc
$\{(x,\psi(x)):x\in\overline{I_\varepsilon}\}$. They satisfy
\begin{equation}\label{eq:strip-flattening}
 c_i(\varepsilon)\longrightarrow u_i(p),\qquad
 \|\zeta_{i,\varepsilon}\|_{H^1(\Omega)}=o(\sqrt\varepsilon),
 \qquad
 \|\zeta_{i,\varepsilon}\|_{L^2(\Omega)}=o(\varepsilon^{3/2}).
\end{equation}
The functions $\zeta_{i,\varepsilon}$ vanish outside
$\Omega\cap B(p,C\varepsilon)$, where $C$ is independent of
$\varepsilon$.
\end{lemma}

\begin{proof}
Fix $p\in \Sigma_*$ and use \eqref{eq:strip-domain}, with $p=0$.
For a sufficiently large fixed number $N_0$, put
\[
 Q_\varepsilon
 =\{(x,y):|x|<N_0\varepsilon,\
            \psi(x)-N_0\varepsilon<y<\psi(x)\},
 \qquad
 c_i(\varepsilon)=\frac1{|Q_\varepsilon|}
                         \int_{Q_\varepsilon}u_i.
\]
The map $(x,t)\mapsto(x,\psi(x)+t)$ has uniformly bounded
bi-Lipschitz constants. The scaled Poincar\'e inequality is therefore
\begin{equation}\label{eq:strip-local-poincare}
 \int_{Q_\varepsilon}|u_i-c_i(\varepsilon)|^2
 \le C\varepsilon^2\int_{Q_\varepsilon}|\nabla u_i|^2.
\end{equation}
Choose $\chi_\varepsilon\in C_c^\infty(\R^2)$ such that
\[
 \begin{gathered}
 0\le\chi_\varepsilon\le1,\qquad
 \chi_\varepsilon=1\quad\text{in a neighborhood of the graph arc
 over }\overline{I_\varepsilon},
 \qquad |\nabla\chi_\varepsilon|\le C\varepsilon^{-1},\\
 \Omega\cap\operatorname{supp}\chi_\varepsilon
       \subset Q_\varepsilon\subset B(p,C\varepsilon).
 \end{gathered}
\]
A radial function about $p$ has these properties after increasing $N_0$
if necessary.
For $\zeta_{i,\varepsilon}
=\chi_\varepsilon(c_i(\varepsilon)-u_i)$, we obtain
\[
 \begin{aligned}
 \|\nabla\zeta_{i,\varepsilon}\|_2^2
 &\le C\int_{Q_\varepsilon}|\nabla u_i|^2
       +C\varepsilon^{-2}
                  \int_{Q_\varepsilon}|u_i-c_i(\varepsilon)|^2\\
 &\le C\int_{\Omega\cap B(p,C\varepsilon)}|\nabla u_i|^2=o(\varepsilon),\\
 \|\zeta_{i,\varepsilon}\|_2^2
 &\le\int_{Q_\varepsilon}|u_i-c_i(\varepsilon)|^2\\
 &\le C\varepsilon^2\int_{\Omega\cap B(p,C\varepsilon)}|\nabla u_i|^2=o(\varepsilon^3).
 \end{aligned}
\]
Finally, $u_i$ is continuous on $\overline\Omega$
\cite[Theorem~3.14(ii)]{Nittka2011}, whence
\[
 |c_i(\varepsilon)-u_i(p)|
 \le\sup_{Q_\varepsilon}|u_i-u_i(p)|\longrightarrow0.
\]
These estimates imply \eqref{eq:strip-flattening}.
\end{proof}

For the remainder of this subsection, fix $\Omega$, an
$L^2(\Omega)$-orthonormal basis $u_1,u_2$ of $E_1(\Omega)$, and
$p\in\Sigma_*$ as in Lemma~\ref{lem:strip-flattening}. Choose $T$ with
\[
 \Lambda_\Omega^{\mathrm N}<T<\Lambda_{\Omega,3}^{\mathrm N}.
\]
Such a number exists because the eigenvalue is double and the spectrum
is discrete. Let $\Omega_\varepsilon$ be the domain
\eqref{eq:strip-domain} at $p$, with a length $h$ satisfying
\eqref{eq:strip-threshold}. The point $p$ and length $h$ remain fixed
as $\varepsilon\downarrow0$. Put
\[
 \mathbf u(p)=\begin{pmatrix}u_1(p)\\u_2(p)\end{pmatrix},
 \qquad k=\sqrt{\Lambda_\Omega^{\mathrm N}},
 \qquad d_h=k\tan(kh).
\]
Corollary~\ref{lem:lipschitz-multiplicity} gives $\mathbf u(p)\ne0$.
Also $kh<\sqrt Th<\pi/2$, so $d_h>0$.

Let $V_\varepsilon$ be the span of the Neumann eigenfunctions on
$\Omega_\varepsilon$ with eigenvalues
$\Lambda_{\Omega_\varepsilon,1}^{\mathrm N}$ and
$\Lambda_{\Omega_\varepsilon,2}^{\mathrm N}$. By
Lemma~\ref{lem:strip-spectrum}, for all sufficiently small
$\varepsilon$,
\[
 0<\Lambda_{\Omega_\varepsilon,1}^{\mathrm N}\le
 \Lambda_{\Omega_\varepsilon,2}^{\mathrm N}<T<
 \Lambda_{\Omega_\varepsilon,3}^{\mathrm N},
 \qquad \dim V_\varepsilon=2.
\]
Set
\[
 B_\varepsilon=a_\varepsilon-\Lambda_\Omega^{\mathrm N}b_\varepsilon,
\]
and define $a_\Omega,b_\Omega,B_\Omega$ analogously. Write
$c_i=c_i(\varepsilon)$ and $\zeta_i=\zeta_{i,\varepsilon}$ for the
constants and corrections in Lemma~\ref{lem:strip-flattening}.
All error terms below may depend on the fixed point $p$ and length $h$.

The flattened eigenfunctions can now be extended through the entrance
and along the strip. The next lemma computes the leading term of their
shifted energy matrix and estimates how closely they satisfy the
eigenvalue equation.

\begin{lemma}[Trial functions and residual estimates]
\label{lem:strip-trial}
There are functions $q_i^\varepsilon\in H^1(\Omega_\varepsilon)$,
$i=1,2$, with $q_i^\varepsilon|_\Omega=u_i+\zeta_i$, such that the
functionals
$r_i^\varepsilon=B_\varepsilon(q_i^\varepsilon,\cdot)$ satisfy
\begin{equation}\label{eq:strip-residual-small}
 \|r_i^\varepsilon\|_{H^1(\Omega_\varepsilon)^*}
                  =o(\sqrt\varepsilon).
\end{equation}
Their shifted energy and mass matrices satisfy
\begin{equation}\label{eq:strip-trial-matrix}
 B_\varepsilon(q_i^\varepsilon,q_j^\varepsilon)
       =-\varepsilon d_hc_ic_j+o(\varepsilon),
 \qquad i,j\in\{1,2\},
\end{equation}
and
\begin{equation}\label{eq:strip-trial-gram}
 b_\varepsilon(q_i^\varepsilon,q_j^\varepsilon)
       =\delta_{ij}+O(\varepsilon).
\end{equation}
\end{lemma}

\begin{proof}
Set
\[
 g(t)=\frac{\cos(k(h-t))}{\cos(kh)},\qquad
 -g''=\Lambda_\Omega^{\mathrm N} g,\quad g(0)=1,\quad g'(h)=0,\quad g'(0)=d_h.
\]
Define
\[
 q_i^\varepsilon=
 \begin{cases}
 u_i+\zeta_i,&\text{on }\Omega,\\
 c_i,&\text{on }A_\varepsilon,\\
 c_i g(t),&\text{at }(x,\beta_\varepsilon+t)\in S_\varepsilon.
 \end{cases}
\]
Matching traces give $q_i^\varepsilon\in H^1(\Omega_\varepsilon)$.
The identity $B_\Omega(u_i,\varphi|_\Omega)=0$ and integration by parts
in the strip yield
\begin{align}\label{eq:strip-residual-identity}
 r_i^\varepsilon(\varphi)
 :=B_\varepsilon(q_i^\varepsilon,\varphi)
 ={}&B_\Omega(\zeta_i,\varphi|_\Omega)
       -\Lambda_\Omega^{\mathrm N} c_i\int_{A_\varepsilon}\varphi
       \nonumber\\
 &-c_i d_h\int_{I_\varepsilon}\varphi(x,\beta_\varepsilon)\,dx.
\end{align}
Since $c_i=O(1)$ and $|A_\varepsilon|=O(\varepsilon^2)$,
\eqref{eq:strip-flattening} and \eqref{eq:strip-base-trace} give
\[
 \begin{aligned}
 |r_i^\varepsilon(\varphi)|
 &\le C\bigl(\|\zeta_i\|_{H^1(\Omega)}
          +|A_\varepsilon|^{1/2}
          +\varepsilon^{3/4}+\varepsilon\bigr)
                 \|\varphi\|_{H^1(\Omega_\varepsilon)}\\
 &=o(\sqrt\varepsilon)\|\varphi\|_{H^1(\Omega_\varepsilon)}.
 \end{aligned}
\]
This proves \eqref{eq:strip-residual-small}.

The eigenfunction equations imply
\[
 B_\Omega(u_i+\zeta_i,u_j+\zeta_j)
       =B_\Omega(\zeta_i,\zeta_j)=o(\varepsilon),
 \qquad
 \int_0^h(g'^2-\Lambda_\Omega^{\mathrm N} g^2)=[gg']_0^h=-d_h.
\]
Adding the contributions of $\Omega$, $A_\varepsilon$, and $S_\varepsilon$ gives
\[
 \begin{aligned}
 B_\varepsilon(q_i^\varepsilon,q_j^\varepsilon)
 &=B_\Omega(\zeta_i,\zeta_j)-\Lambda_\Omega^{\mathrm N} c_ic_j|A_\varepsilon|
       +\varepsilon c_ic_j\int_0^h(g'^2-\Lambda_\Omega^{\mathrm N} g^2)\\
 &=-\varepsilon d_hc_ic_j+o(\varepsilon).
 \end{aligned}
\]
Similarly, \eqref{eq:strip-flattening} gives
\[
 \begin{aligned}
 b_\varepsilon(q_i^\varepsilon,q_j^\varepsilon)
 &=\delta_{ij}+o(\varepsilon^{3/2})
      +c_ic_j\left(|A_\varepsilon|+
                      \varepsilon\int_0^h g^2\right)\\
 &=\delta_{ij}+O(\varepsilon).
 \end{aligned}
\]

These are \eqref{eq:strip-trial-matrix} and
\eqref{eq:strip-trial-gram}.
\end{proof}

The residual estimate and the spectral gap control the projection of
the trial functions onto $V_\varepsilon$. The resulting error in the
shifted energy is quadratic, so the first-order term survives this
projection and determines the matrix in the next lemma.

\begin{lemma}[Effective matrix]
\label{lem:strip-matrix}
For all sufficiently small $\varepsilon$, there is an
$L^2(\Omega_\varepsilon)$-orthonormal basis
$(f_1^\varepsilon,f_2^\varepsilon)$ of $V_\varepsilon$ such that
$f_i^\varepsilon|_\Omega\to u_i$ in $H^1(\Omega)$. The matrix
$H_\varepsilon=
\bigl(a_\varepsilon(f_i^\varepsilon,f_j^\varepsilon)\bigr)_{i,j=1}^2$
satisfies
\begin{equation}\label{eq:strip-effective-matrix}
 H_\varepsilon=\Lambda_\Omega^{\mathrm N}I
       -\varepsilon d_h\mathbf u(p)\mathbf u(p)^T+o(\varepsilon).
\end{equation}
\end{lemma}

\begin{proof}
Let $q_i^\varepsilon$ and $r_i^\varepsilon$ be as in
Lemma~\ref{lem:strip-trial}. Let $P_\varepsilon$ be the
$L^2(\Omega_\varepsilon)$-orthogonal projection onto $V_\varepsilon$,
and choose an
$L^2$-orthonormal Neumann eigenbasis
$\{\phi_m^\varepsilon\}_{m\ge0}$ with eigenvalues
$\Lambda_{\Omega_\varepsilon,m}^{\mathrm N}$. For
$\omega=\sum_m\omega_m\phi_m^\varepsilon\in H^1(\Omega_\varepsilon)$
and the functional $\ell_\omega=B_\varepsilon(\omega,\cdot)$,
\[
 \|\omega\|_{H^1}^2=\sum_m(1+\Lambda_{\Omega_\varepsilon,m}^{\mathrm N})|\omega_m|^2,
 \qquad
 \|\ell_\omega\|_{H^1(\Omega_\varepsilon)^*}^2
 =\sum_m\frac{|\Lambda_{\Omega_\varepsilon,m}^{\mathrm N}-\Lambda_\Omega^{\mathrm N}|^2}
                 {1+\Lambda_{\Omega_\varepsilon,m}^{\mathrm N}}|\omega_m|^2.
\]
By Lemma~\ref{lem:strip-spectrum}, there are $\varepsilon_0>0$ and
$C_*<\infty$ such that
\[
 \sup_{\substack{m\notin\{1,2\}\\0<\varepsilon<\varepsilon_0}}
 \frac{1+\Lambda_{\Omega_\varepsilon,m}^{\mathrm N}}{|\Lambda_{\Omega_\varepsilon,m}^{\mathrm N}-\Lambda_\Omega^{\mathrm N}|}\le C_*.
\]
Indeed, for $m=0$ the quotient equals $1/\Lambda_\Omega^{\mathrm N}$,
and for $m\ge3$ it is at most
$(1+T)/(T-\Lambda_\Omega^{\mathrm N})$, because
$\Lambda_{\Omega_\varepsilon,m}^{\mathrm N}>T$. Put
$e_i^\varepsilon=(I-P_\varepsilon)q_i^\varepsilon$. The coefficients of
$e_i^\varepsilon$ vanish for $m\in\{1,2\}$ and equal those of
$q_i^\varepsilon$ for $m\notin\{1,2\}$, and for such $m$ we have
$1+\Lambda_{\Omega_\varepsilon,m}^{\mathrm N}\le
C_*^2|\Lambda_{\Omega_\varepsilon,m}^{\mathrm N}-\Lambda_\Omega^{\mathrm N}|^2/
(1+\Lambda_{\Omega_\varepsilon,m}^{\mathrm N})$. Consequently, since
$r_i^\varepsilon=\ell_\omega$ for $\omega=q_i^\varepsilon$,
\begin{equation}\label{eq:strip-projection-error}
 \|e_i^\varepsilon\|_{H^1(\Omega_\varepsilon)}^2
 \le C_*^2\|r_i^\varepsilon\|_{H^1(\Omega_\varepsilon)^*}^2
 =o(\varepsilon).
\end{equation}
This squared projection error is smaller than the order-$\varepsilon$
term in the strip expansion.
Put $z_i^\varepsilon=P_\varepsilon q_i^\varepsilon$. Orthogonality
for both $a_\varepsilon$ and $b_\varepsilon$ yields
\[
 \begin{aligned}
 B_\varepsilon(z_i^\varepsilon,z_j^\varepsilon)
 &=B_\varepsilon(q_i^\varepsilon,q_j^\varepsilon)
       -B_\varepsilon(e_i^\varepsilon,e_j^\varepsilon)
 =-\varepsilon d_hc_ic_j+o(\varepsilon),\\
 G_\varepsilon
 &:=(b_\varepsilon(z_i^\varepsilon,z_j^\varepsilon))_{i,j=1}^2
   =I+O(\varepsilon).
 \end{aligned}
\]
Thus $G_\varepsilon>0$ for small $\varepsilon$, and
\[
 (f_1^\varepsilon,f_2^\varepsilon)
    =(z_1^\varepsilon,z_2^\varepsilon)G_\varepsilon^{-1/2}
\]
is an $L^2(\Omega_\varepsilon)$-orthonormal basis of $V_\varepsilon$.
Since $a_\varepsilon=B_\varepsilon+\Lambda_\Omega^{\mathrm N}b_\varepsilon$,
the matrix
$H_\varepsilon=\bigl(a_\varepsilon(f_i^\varepsilon,f_j^\varepsilon)\bigr)$
is
\[
 \begin{aligned}
 H_\varepsilon
 &=\Lambda_\Omega^{\mathrm N} I+G_\varepsilon^{-1/2}
    (B_\varepsilon(z_i^\varepsilon,z_j^\varepsilon))_{i,j=1}^2
                         G_\varepsilon^{-1/2}\\
 &=\Lambda_\Omega^{\mathrm N} I-\varepsilon d_h\mathbf u(p)\mathbf u(p)^T+o(\varepsilon),
 \end{aligned}
\]
because $(c_1,c_2)^T\to \mathbf u(p)$.

Since
$z_i^\varepsilon|_\Omega=u_i+\zeta_i-e_i^\varepsilon|_\Omega$ and
$G_\varepsilon\to I$, \eqref{eq:strip-flattening} and
\eqref{eq:strip-projection-error} also give
\[
 f_i^\varepsilon|_\Omega\longrightarrow u_i
       \quad\text{in }H^1(\Omega).
\]
\end{proof}

The leading perturbation matrix has one negative eigenvalue and one
zero eigenvalue. It therefore identifies which linear combination
becomes the first eigenfunction when the double eigenvalue splits.

\begin{proposition}[Splitting and selection of the first eigenfunction]
\label{prop:strip-evaluation}
In the fixed setting above,
\begin{equation}\label{eq:strip-eigenvalue-splitting}
 \Lambda_{\Omega_\varepsilon,1}^{\mathrm N}
       =\Lambda_\Omega^{\mathrm N}
        -\varepsilon d_h|\mathbf u(p)|^2+o(\varepsilon),
 \qquad
 \Lambda_{\Omega_\varepsilon,2}^{\mathrm N}
       =\Lambda_\Omega^{\mathrm N}+o(\varepsilon).
\end{equation}
For all sufficiently small $\varepsilon$, the first positive
Neumann eigenvalue of $\Omega_\varepsilon$ is simple. Its normalized
eigenfunction can be chosen so that its restriction to $\Omega$
converges in $H^1(\Omega)$ and in $C^\infty_{\mathrm{loc}}(\Omega)$ to
\begin{equation}\label{eq:strip-evaluation-direction}
 w_p=\frac{u_1(p)u_1+u_2(p)u_2}{|\mathbf u(p)|}.
\end{equation}
\end{proposition}

\begin{proof}
Use the basis and matrix supplied by Lemma~\ref{lem:strip-matrix}.
The eigenvalues of $H_\varepsilon$
are $\Lambda_{\Omega_\varepsilon,1}^{\mathrm N}$ and
$\Lambda_{\Omega_\varepsilon,2}^{\mathrm N}$, because $H_\varepsilon$ is
the matrix of $a_\varepsilon$ on $V_\varepsilon$ in an orthonormal
basis. We have
\[
 \varepsilon^{-1}(H_\varepsilon-\Lambda_\Omega^{\mathrm N} I)
        \longrightarrow-d_h\mathbf u(p)\mathbf u(p)^T,
 \qquad
 \operatorname{spec}(-d_h\mathbf u(p)\mathbf u(p)^T)=\{-d_h|\mathbf u(p)|^2,0\}.
\]
The limiting eigenvalues are distinct, because $\mathbf u(p)\ne0$ and
$d_h>0$. Hence \eqref{eq:strip-eigenvalue-splitting} holds, and
$\Lambda_{\Omega_\varepsilon,1}^{\mathrm N}<\Lambda_{\Omega_\varepsilon,2}^{\mathrm N}$
for small $\varepsilon$, so the first positive eigenvalue of
$\Omega_\varepsilon$ is simple. A unit eigenvector
$\xi_\varepsilon\in\R^2$ of $H_\varepsilon$ for its smaller eigenvalue
satisfies, after a sign choice,
\[
 \xi_\varepsilon\longrightarrow \mathbf u(p)/|\mathbf u(p)|.
\]
By the $H^1(\Omega)$ convergence in
Lemma~\ref{lem:strip-matrix},
\[
 \sum_{i=1}^2\xi_{\varepsilon,i}f_i^\varepsilon|_\Omega
       \longrightarrow
       \frac{u_1(p)u_1+u_2(p)u_2}{|\mathbf u(p)|}=w_p
       \quad\text{in }H^1(\Omega).
\]
The function $\sum_i\xi_{\varepsilon,i}f_i^\varepsilon$ is a normalized
first eigenfunction on $\Omega_\varepsilon$.
The last convergence also holds in $C^\infty_{\mathrm{loc}}(\Omega)$
by interior elliptic estimates and
$\Lambda_{\Omega_\varepsilon,1}^{\mathrm N}\to\Lambda_\Omega^{\mathrm N}$.
\end{proof}

\subsection{Approximation of a prescribed eigenfunction}
\label{subsec:prescribed-eigenfunction}

The selected direction is determined by boundary evaluation at the
attachment point. In the double-eigenvalue case, a boundary zero of
an eigenfunction orthogonal to the prescribed $u$ gives the desired
direction, because $u$ does not vanish there. Density of admissible
points allows us to approach that boundary point.

\begin{proposition}[Approximation of a prescribed eigenfunction]
\label{prop:prescribed}
Let $\Omega\subset\R^2$ be a bounded simply connected Lipschitz domain.
For every normalized $u\in E_1(\Omega)$, there are
bounded simply connected Lipschitz domains $\Omega_j\supset\Omega$ with
simple first positive Neumann eigenvalue and normalized
eigenfunctions $u_j\in E_1(\Omega_j)$ such that
\begin{equation}\label{eq:lipschitz-selected-convergence}
 \Lambda_{\Omega_j}^{\mathrm N}\longrightarrow\Lambda_\Omega^{\mathrm N},\qquad
 u_j|_\Omega\longrightarrow u\text{ in }H^1(\Omega)
       \text{ and in }C^\infty_{\mathrm{loc}}(\Omega).
\end{equation}

\end{proposition}

\begin{proof}
If $\dim E_1(\Omega)=1$, take $(\Omega_j,u_j)=(\Omega,u)$. Otherwise,
Corollary~\ref{lem:lipschitz-multiplicity} gives $\dim E_1(\Omega)=2$, so
\[
 E_1(\Omega)=\operatorname{span}\{u,v\},\qquad
 \langle u,v\rangle_{L^2(\Omega)}=0,\qquad \|v\|_{L^2(\Omega)}=1.
\]
By Lemma~\ref{lem:nodal-arc}, the nodal arc of $v$ has an endpoint
$p_0\in\partial\Omega$. Corollary~\ref{lem:lipschitz-multiplicity}
shows that the basis $u,v$ has no common zero on $\partial\Omega$, so
\[
 v(p_0)=0,\qquad u(p_0)\ne0.
\]
Choose $p_j\in \Sigma_*$ with $p_j\to p_0$. This is possible because
$\Sigma_*$ is dense (Lemma~\ref{lem:strip-good-points}). Set
$s=\operatorname{sgn}u(p_0)$. Continuity of $u$ and $v$ on
$\overline\Omega$ yields
\begin{equation}\label{eq:strip-direction-density}
 w_j=s\,
 \frac{u(p_j)u+v(p_j)v}{\sqrt{u(p_j)^2+v(p_j)^2}}
 \longrightarrow u\qquad\text{in }H^1(\Omega).
\end{equation}
The denominator is nonzero by Corollary~\ref{lem:lipschitz-multiplicity}.

Since the first positive eigenvalue is double and the spectrum is
discrete, we may fix
$\Lambda_\Omega^{\mathrm N}<T<\Lambda_{\Omega,3}^{\mathrm N}$.
For each $j$, choose
$h_j>0$ so small that the construction \eqref{eq:strip-domain} at $p_j$
is defined and $\pi^2/(4h_j^2)>T$. Apply
Proposition~\ref{prop:strip-evaluation} with $(u_1,u_2)=(u,v)$,
$p=p_j$, and $h=h_j$, and multiply the eigenfunction that it supplies by
$s$. This gives a width $\varepsilon_j>0$ such that the domain
$\Omega_j$, obtained by attaching to $\Omega$ the strip of width
$\varepsilon_j$ and length $h_j$ at $p_j$, and a normalized
eigenfunction $u_j\in E_1(\Omega_j)$ satisfy
\[
 \Lambda_{\Omega_j,1}^{\mathrm N}<\Lambda_{\Omega_j,2}^{\mathrm N},\qquad
 |\Lambda_{\Omega_j}^{\mathrm N}-\Lambda_\Omega^{\mathrm N}|<j^{-1},\qquad
 \|u_j|_\Omega-w_j\|_{H^1(\Omega)}<j^{-1}.
\]
Here $p_j,h_j$ are fixed before $\varepsilon_j$ is chosen. Consequently,
\[
 \|u_j|_\Omega-u\|_{H^1(\Omega)}
 \le j^{-1}+\|w_j-u\|_{H^1(\Omega)}\longrightarrow0.
\]
Since $-\Delta u_j=\Lambda_{\Omega_j}^{\mathrm N}u_j$ on $\Omega$ and
$\Lambda_{\Omega_j}^{\mathrm N}\to\Lambda_\Omega^{\mathrm N}$, interior elliptic estimates give
$u_j\to u$ in $C^\infty_{\mathrm{loc}}(\Omega)$.
\end{proof}

\section{Smooth approximation and strict boundary attainment}
\label{sec:strict-attainment}

We now pass from the domains of Proposition~\ref{prop:prescribed} to
smooth domains. Theorem~\ref{thm:smooth} is invoked only after this
approximation has been established.

\subsection{Conformal exhaustion and eventual simplicity}
\label{subsec:exhaustion}

We first put the Neumann problems on a fixed disk, then prove spectral
convergence, and finally select the eigenfunction when the limiting
eigenvalue is simple.

Throughout this subsection, let $\Omega$ be a bounded simply connected
Lipschitz domain and let $F:\D\to\Omega$ be conformal. For $0<r<1$, put
\[
 F_r(z)=F(rz),\qquad G_r=F(r\D),\qquad
 \rho=|F'|^2,\qquad \rho_r=|F_r'|^2=r^2|F'(rz)|^2.
\]
Since $F$ is univalent and holomorphic on a neighborhood of
$\overline{r\D}$, with $F'\ne0$, the domains $G_r$ are simply connected
and $C^\infty$. They satisfy $\overline{G_r}\subset\Omega$ and
$G_r\uparrow\Omega$ as $r\uparrow1$.

For $f,g\in H^1(\D)$, write
\[
 e(f,g)=\int_\D\nabla f\cdot\nabla g,\qquad
 b(f,g)=\int_\D fg\rho,\qquad
 b_r(f,g)=\int_\D fg\rho_r,\qquad
 \langle f,g\rangle_*=e(f,g)+b(f,g).
\]

\begin{lemma}[Conformal pullback and convergence of the mass forms]
\label{lem:conformal-forms}
The maps $u\mapsto u\circ F$ and $u\mapsto u\circ F_r$ are bijections
from $H^1(\Omega)$ and $H^1(G_r)$, respectively, onto $H^1(\D)$.
They transform the Neumann energy and mass forms into $(e,b)$ and
$(e,b_r)$. The inner product $\langle\cdot,\cdot\rangle_*$ is equivalent
to the usual $H^1(\D)$ inner product. For some $p>1$,
\begin{equation}\label{eq:bridge-weight-convergence}
 \rho_r\longrightarrow\rho\quad\text{in }L^p(\D),
\end{equation}
and the forms $b_r$ converge to $b$ in operator norm on $H^1(\D)$.
\end{lemma}

\begin{proof}
By \eqref{eq:lipschitz-conformal-integrability},
$\rho\in L^p(\D)$ for some $p>1$. Continuity of dilations in
$L^{2p}(\D)$ gives \eqref{eq:bridge-weight-convergence}.
With $p'=p/(p-1)$, H\"older's inequality and the Sobolev embedding give
\begin{equation}\label{eq:bridge-mass-form-convergence}
 |b_r(f,g)-b(f,g)|
 \le\|\rho_r-\rho\|_p\|f\|_{2p'}\|g\|_{2p'}
 \le C\|\rho_r-\rho\|_p\|f\|_{H^1}\|g\|_{H^1}.
\end{equation}
The same estimate bounds $e+b$ from above by the $H^1$ norm squared.
The lower bound follows from the anchored Poincar\'e inequality
\[
 \|f\|_{L^2(\D)}^2
 \le C\left(\|\nabla f\|_{L^2(\D)}^2+
                   \|f\|_{L^2(\frac12\D)}^2\right),
 \qquad \min_{\overline{\frac12\D}}\rho>0.
\]

If $f\in H^1(\D)$, conformal invariance shows that $f\circ F^{-1}$
has Dirichlet energy $e(f,f)$ and squared $L^2(\Omega)$ norm
$b(f,f)<\infty$. Hence $f\circ F^{-1}\in H^1(\Omega)$.
Conversely, for $u\in H^1(\Omega)$, the function $u\circ F$ is locally
in $H^1$, has finite energy, and is square integrable on $\frac12\D$.
Apply the anchored Poincar\'e inequality on $r\D$. By scaling, its
constant is independent of $r\ge3/4$. Monotone convergence then gives
$u\circ F\in H^1(\D)$. The same argument applies to $F_r$ and $G_r$.

The changes of variables preserve the energy and give the stated
weighted mass forms. Thus $\gamma$ is a Neumann eigenvalue of $\Omega$
(respectively of $G_r$) if and only if, for some $f\ne0$,
\[
 e(f,g)=\gamma b(f,g)
 \quad\text{(respectively }e(f,g)=\gamma b_r(f,g)\text{)}
 \qquad(g\in H^1(\D)).
\]
\end{proof}

We now express the convergence of the mass forms as convergence of
compact self-adjoint operators on a fixed Hilbert space.
Equip $H^1(\D)$ with the inner product $\langle\cdot,\cdot\rangle_*$
and denote this Hilbert space by $\mathcal H$. Define $B$ and $B_r$ by
\[
 \langle Bf,g\rangle_*=b(f,g),\qquad
 \langle B_rf,g\rangle_*=b_r(f,g),
\]
and set
\[
 A_r=I+B_r-B,\qquad C_r=A_r^{-1/2}B_rA_r^{-1/2}
\]
whenever $A_r$ is positive and invertible.

\begin{lemma}[Spectral convergence along conformal exhaustion]
\label{lem:conformal-spectrum}
For $r$ sufficiently close to one, $A_r$ is positive and invertible,
and $C_r\to B$ in operator norm on $\mathcal H$. Both $C_r$ and $B$
are compact and self-adjoint. Their nonzero eigenvalues, in decreasing
order and with multiplicities, are
\[
 \mu_k(C_r)=(1+\Lambda_{G_r,k}^{\mathrm N})^{-1},\qquad
 \mu_k(B)=(1+\Lambda_{\Omega,k}^{\mathrm N})^{-1},\qquad k\ge0.
\]
Consequently $\Lambda_{G_r,k}^{\mathrm N}\to\Lambda_{\Omega,k}^{\mathrm N}$
for every $k\ge0$.
\end{lemma}

\begin{proof}
Choose $p>1$ as in Lemma~\ref{lem:conformal-forms} and put
$p'=p/(p-1)$. The compact embedding
$H^1(\D)\hookrightarrow L^{2p'}(\D)$ shows that
$B$ and $B_r$ are compact, nonnegative, and self-adjoint. By
\eqref{eq:bridge-mass-form-convergence}, $\|B_r-B\|\to0$.
It follows that $A_r\to I$, so $A_r$ is positive and invertible for
$r$ close to one, and $C_r\to B$ in norm. Also,
\[
 \langle A_rf,g\rangle_*=e(f,g)+b_r(f,g).
\]
Under the change of variable $z=A_r^{1/2}f$, we have
\[
 e(f,g)=\gamma b_r(f,g)\quad(g\in H^1(\D))
 \quad\Longleftrightarrow\quad
 C_rz=(1+\gamma)^{-1}z.
\]
The same computation for $b$ replaces $A_r$ by $I$ and $C_r$ by $B$.
Lemma~\ref{lem:conformal-forms} therefore identifies their nonzero
eigenvalues as stated. Norm convergence of compact self-adjoint
operators gives convergence of each ordered eigenvalue, and hence
of the corresponding Neumann eigenvalues.
\end{proof}

The domains supplied by Proposition~\ref{prop:prescribed} have a simple
first positive eigenvalue. The next lemma uses the corresponding
spectral projection to approximate the prescribed eigenfunction on
the smooth exhaustion.

\begin{lemma}[Eventual simplicity along conformal exhaustion]
\label{lem:eventual-simplicity}
Suppose that $\Lambda_\Omega^{\mathrm N}$ is simple. There is $r_*<1$
such that
\begin{equation}\label{eq:eventual-simple}
 \Lambda_{G_r}^{\mathrm N}\longrightarrow\Lambda_\Omega^{\mathrm N}
 \quad(r\uparrow1),
 \qquad \dim E_1(G_r)=1\quad\text{for }r_*<r<1.
\end{equation}
For every normalized $v\in E_1(\Omega)$ there are normalized
$v_r\in E_1(G_r)$, $r_*<r<1$, such that
\begin{equation}\label{eq:smooth-local-convergence}
 v_r\longrightarrow v\quad\text{in }C^\infty_{\mathrm{loc}}(\Omega)
 \qquad(r\uparrow1).
\end{equation}
\end{lemma}

\begin{proof}
By Lemma~\ref{lem:conformal-spectrum},
$\mu_*=(1+\Lambda_\Omega^{\mathrm N})^{-1}$ is a simple isolated
eigenvalue of $B$. Let $P$ be its spectral projection, and let $P_r$
be the projection of $C_r$ associated with a fixed small disk about
$\mu_*$ containing no other spectral point of $B$.
Norm convergence gives
\[
 \|P_r-P\|\longrightarrow0,\qquad \operatorname{rank}P=1.
\]
The same lemma gives $\Lambda_{G_r}^{\mathrm N}\to\Lambda_\Omega^{\mathrm N}$.
We may therefore choose $r_*<1$ such that, for $r_*<r<1$,
$\|P_r-P\|<1$ and $(1+\Lambda_{G_r}^{\mathrm N})^{-1}$ lies in the
disk defining $P_r$. Then $\operatorname{rank}P_r=1$, which proves
eventual simplicity.

For the prescribed normalized $v$, put $w=v\circ F$. Then
$Pw=w$, $b(w,w)=1$, and
\[
 \|P_rw-w\|_*<\|w\|_*,\qquad P_rw\ne0.
\]
Define
\[
 \widetilde w_r=A_r^{-1/2}P_rw,\qquad
 w_r=\frac{\widetilde w_r}
 {b_r(\widetilde w_r,\widetilde w_r)^{1/2}}.
\]
The denominator tends to one, and $w_r\to w$ in $H^1(\D)$.
Moreover,
\[
 -\Delta w_r=\Lambda_{G_r}^{\mathrm N}\rho_rw_r,\qquad
 \rho_r\longrightarrow\rho\quad\text{in }C^\infty_{\mathrm{loc}}(\D).
\]
Interior elliptic estimates give $w_r\to w$ in
$C^\infty_{\mathrm{loc}}(\D)$. The normalized eigenfunctions
\[
 v_r(y)=w_r(F^{-1}(y)/r),\qquad y\in G_r,
\]
belong to $E_1(G_r)$ and converge to $v$ in
$C^\infty_{\mathrm{loc}}(\Omega)$.
\end{proof}

The exponent $p$ and the constants above depend on $\Omega$.
This is sufficient, because the proof of Theorem~\ref{thm:lipschitz}
exhausts one domain $\Omega_j$ at a time.

\subsection{Proof of Theorem~\ref{thm:lipschitz}}
\label{subsec:lipschitz-proof}

The approximation will exclude strict interior local extrema.
To obtain (HS1), we must also exclude global extrema attained along
curves. The next lemma supplies this step through the critical-set
argument of Judge and Mondal
\cite[Proposition~14, p.~160]{JudgeMondalCritical2025}.
That proposition excludes an interior critical hypersurface for a
first positive Neumann eigenfunction on a bounded simply connected
Lipschitz domain, without a simplicity assumption. We give the planar
argument needed here, using the weak eigenvalue equation in place of
integration by parts on a component cut out by a critical arc.

\begin{lemma}[Strict boundary attainment]
\label{lem:strict-boundary-attainment}
Let $\Omega\subset\R^2$ be a bounded simply connected Lipschitz domain,
and let $0\ne u\in E_1(\Omega)$. If $u$ has no strict local maximum
or minimum in $\Omega$, then
\[
 \min_{\partial\Omega}u<u(x)<\max_{\partial\Omega}u
 \qquad(x\in\Omega).
\]
\end{lemma}

\begin{proof}
Suppose that $u(x_0)=M:=\max_{\overline\Omega}u$ for some
$x_0\in\Omega$, and put
\[
 Z=\{x\in\Omega:u(x)=M\}.
\]
Since $u\ne0$ and $\int_\Omega u=0$, we have $M>0$. At every $z\in Z$,
\[
 \nabla u(z)=0,
 \qquad \Delta u(z)=-\Lambda_\Omega^{\mathrm N}M<0.
\]
Thus $\partial_{jj}u(z)\ne0$ for some $j\in\{1,2\}$. Interior
analyticity and the analytic implicit-function theorem show that
$\{\partial_j u=0\}$ is a real-analytic arc near $z$, containing $Z$
there. The restriction of $u-M$ to this arc is real-analytic. Either
$z$ is an isolated point of $Z$, or that restriction vanishes
identically near $z$. The first alternative would make $z$ a strict
local maximum. Hence $Z$ is a relatively closed, embedded
real-analytic one-dimensional manifold without boundary in $\Omega$.

Choose a connected component $K$ of $Z$. It is closed in $\Omega$,
so its inclusion is proper. Since $\Omega$ is simply connected,
$\Omega\setminus K$ has at least two connected components by the
separation theorem for proper embedded hypersurfaces
\cite[Lemma~7, p.~158]{JudgeMondalCritical2025}.

Let $U$ be any such component and define
\[
 \phi=(M-u)\mathbf1_U.
\]
We claim that
\begin{equation}\label{eq:critical-arc-test}
 \phi\in H^1(\Omega),
 \qquad \nabla\phi=-\mathbf1_U\nabla u.
\end{equation}
Indeed, $\partial U\cap\Omega\subset K$, and on $K$ both $M-u$ and
$\nabla u$ vanish. For $z\in K$,
\[
 |\phi(x)-\phi(z)|\le |M-u(x)|=o(|x-z|).
\]
Off $K$, the indicator of $U$ is locally constant. It follows that
$\phi\in C^1(\Omega)$ and that its gradient is as in
\eqref{eq:critical-arc-test}, with value zero on $K$. Moreover,
\[
 \|\phi\|_{L^2(\Omega)}
 \le M|\Omega|^{1/2}+\|u\|_{L^2(\Omega)},
 \qquad
 \|\nabla\phi\|_{L^2(\Omega)}
 \le\|\nabla u\|_{L^2(\Omega)}.
\]
This proves the claim without any regularity assumption on
$\partial U\cap\partial\Omega$.

Testing the weak Neumann equation with $\phi$ gives
\begin{equation}\label{eq:critical-arc-sign}
 -\int_U|\nabla u|^2\,dx
 =\Lambda_\Omega^{\mathrm N}\int_Uu(M-u)\,dx<0.
\end{equation}
The strict inequality holds because vanishing of the gradient on the
nonempty open set $U$ would, by the eigenvalue equation and unique
continuation, imply $u\equiv0$. Since $M-u\ge0$, equation
\eqref{eq:critical-arc-sign} shows that $U$ contains a point where
$u<0$.

Thus $\{u<0\}$ meets two different components of $\Omega\setminus K$.
It cannot connect them, because $u=M>0$ on $K$. This contradicts the
connectedness of the negative nodal domain in
Lemma~\ref{lem:two-nodal-domains}. Therefore the global
maximum is attained only on $\partial\Omega$. Applying the argument
to $-u$ proves the assertion for the minimum.
\end{proof}

We can now combine the two approximations. Proposition~\ref{prop:prescribed}
selects the given eigenfunction, and Lemma~\ref{lem:eventual-simplicity}
approximates it on smooth domains. The smooth theorem rules out compact
interior superlevel components, after which
Lemma~\ref{lem:strict-boundary-attainment} gives (HS1).

\begin{proof}[Proof of Theorem~\ref{thm:lipschitz}]
Let $0\ne u\in E_1(\Omega)$. After normalization we may suppose
$\|u\|_{L^2(\Omega)}=1$. Suppose that, for some $t\in\R$, a connected
component $U$ of $\{x\in\Omega:u(x)>t\}$ has compact closure in
$\Omega$. Fix $x_0\in U$ and put
\[
 \eta=u(x_0)-t>0.
\]
Together with the compactness assumption, continuity gives
\begin{equation}\label{eq:compact-superlevel}
 \overline U\text{ is a compact subset of }\Omega,
 \qquad u=t\text{ on }\partial U.
\end{equation}
Indeed, at any $y\in\partial U$, continuity gives
$u(y)\ge t$. If $u(y)>t$, a ball about $y$ is contained in
$\{u>t\}$ and intersects $U$, so it lies in the same component,
contrary to $y\in\partial U$. This proves
\eqref{eq:compact-superlevel}. No regularity of $\partial U$ is used.

Let $(\Omega_j,u_j)$ be the sequence supplied by
Proposition~\ref{prop:prescribed}. Since $u_j\to u$ in
$C^\infty_{\mathrm{loc}}(\Omega)$, we can choose $j$ so large that
\[
 \|u_j-u\|_{C(\overline U)}<\frac\eta4.
\]
The first positive eigenvalue of $\Omega_j$ is simple.
Lemma~\ref{lem:eventual-simplicity}, applied to $\Omega_j$ with a
conformal map $F=F_j:\D\to\Omega_j$ and $v=u_j$, gives $r_*<1$, smooth
simply connected domains $G_{j,r}=F_j(r\D)$ with simple first positive
eigenvalue, and normalized eigenfunctions $v_{j,r}\in E_1(G_{j,r})$,
$r_*<r<1$, such that $v_{j,r}\to u_j$ in
$C^\infty_{\mathrm{loc}}(\Omega_j)$ as $r\uparrow1$.
Choose $r$ such that
\[
 \overline U\subset G_{j,r},
 \qquad\|v_{j,r}-u_j\|_{C(\overline U)}<\frac\eta4.
\]
The two error bounds and \eqref{eq:compact-superlevel} imply
\begin{equation}\label{eq:first-superlevel-gap}
 \max_{\partial U}v_{j,r}<t+\frac\eta2<v_{j,r}(x_0).
\end{equation}
The continuous function $v_{j,r}$ therefore attains its maximum
on the compact set $\overline U$ at a point of $U$. That point is
an interior local maximum of $v_{j,r}$ in $G_{j,r}$, so the gradient
of $v_{j,r}$ vanishes there. This contradicts
Theorem~\ref{thm:smooth}. Thus no component of $\{u>t\}$ has compact
closure in $\Omega$, for any $t\in\R$. Applying the same argument to
$-u$ also excludes compactly contained components of $\{u<t\}$.

Consequently $u$ has no strict interior local extremum. Indeed, at a
strict local maximum $x_0$ we could choose
$\overline{B(x_0,r)}\subset\Omega$ and $t$ such that
\[
 \max_{\partial B(x_0,r)}u<t<u(x_0).
\]
The component of $\{u>t\}$ containing $x_0$ would then have compact
closure in $B(x_0,r)$, a contradiction. The same reasoning applies
to a strict local minimum. Lemma~\ref{lem:strict-boundary-attainment}
now gives \eqref{eq:HS1} and, in particular, both identities in
\eqref{eq:main-HS2}.
\end{proof}

\begin{remark}[Compact level-set components]
\label{rem:superlevel}
The proof also shows that, for every $u\in E_1(\Omega)$ and
$t\in\R$, no connected component of $\{u>t\}$ or $\{u<t\}$ has compact
closure in $\Omega$. This is the approximation input used in
Lemma~\ref{lem:strict-boundary-attainment} through the exclusion of
strict local extrema. Property (HS2) alone does not provide that
exclusion. The conclusion (HS1) rules out interior global extrema;
it does not assert that $\nabla u$ is everywhere nonzero in $\Omega$.
\end{remark}

\section{Polygonal domains and starlike examples}
\label{sec:polygonal-comments}

\subsection{Polygonal domains}
\label{subsec:polygons}

For bounded simply connected polygonal domains, known nodal results
give the nodal-arc conclusion of Lemma~\ref{lem:nodal-arc} and the
multiplicity bound of Corollary~\ref{lem:lipschitz-multiplicity}.
Corner expansions show that every boundary point is admissible.
Thus the approximation of the point $p_0$ in the proof of
Proposition~\ref{prop:prescribed} is not needed.

Let $\Omega$ be a bounded simply connected polygonal domain. At a
vertex $z_0$ of interior angle $\beta\in(0,2\pi)$, put
$\alpha=\pi/\beta>1/2$. In polar coordinates centered at $z_0$, with
$\theta$ measured from one of the two sides at $z_0$, every
$u\in E_1(\Omega)$ has the expansion
\begin{equation}\label{eq:polygon-corner-series}
 u(r,\theta)=\sum_{m\ge0}c_m
 J_{m\alpha}\!\left(\sqrt{\Lambda_\Omega^{\mathrm N}}\,r\right)
 \cos(m\alpha\theta),\qquad 0\le\theta\le\beta,
\end{equation}
near $z_0$, as in
\cite[Section~4, pp.~178--179, equations~(5)--(6)]{JudgeMondal}. We
call $c_0$ and $c_1$ the first two \emph{corner coefficients} of $u$ at
$z_0$. In particular, $c_0=u(z_0)$. A point in the interior of a side is
regarded as a vertex with $\beta=\pi$.

Judge and Mondal
\cite[Theorem~5.2, p.~182, and Corollaries~5.3--5.4, p.~183]{JudgeMondal}
prove, without a convexity assumption, the description of the nodal set
in Lemma~\ref{lem:nodal-arc} and the bound
$\dim E_1(\Omega)\le2$. Their Corollary~5.3 also shows that
$c_0=c_1=0$ implies $u=0$. Thus the map $u\mapsto(c_0,c_1)$ is
injective on $E_1(\Omega)$ at each vertex $z_0$.
Thus the reflection to the sphere and Kokarev's theorem
\cite{Kokarev2014} are not needed for polygons.

For a polygon, every boundary point is admissible. Indeed,
differentiating the expansion \eqref{eq:polygon-corner-series} and
using $J_0'(r)=O(r)$ gives, at a vertex $z_0$,
\begin{equation}\label{eq:polygon-corner-energy}
 |\nabla u|=O(r^{\alpha-1}+r),\qquad
 \int_{\Omega\cap B(z_0,r)}|\nabla u|^2
       =O(r^{2\alpha}+r^4)=o(r).
\end{equation}
At a point in the interior of a side, reflection gives bounded
$\nabla u$ and the stronger estimate $O(r^2)$. Hence
\eqref{eq:strip-good-energy} holds at \emph{every} boundary point,
including reentrant vertices; that is, $\Sigma_*=\partial\Omega$.

Suppose that $\dim E_1(\Omega)=2$. Choose an orthonormal basis $u,v$ of
$E_1(\Omega)$ and an endpoint $p_0$ of the nodal arc of $v$. Then
$u(p_0)\ne0$: the first corner coefficient at $p_0$ is the value at
$p_0$, so otherwise it would vanish on all of $E_1(\Omega)$, which
contradicts the injectivity of the map $u\mapsto(c_0,c_1)$ on the
two-dimensional space $E_1(\Omega)$. We may therefore take $p=p_0$ in
Proposition~\ref{prop:strip-evaluation}. The limit in that proposition
is
\[
 \frac{u(p_0)u+v(p_0)v}{\sqrt{u(p_0)^2+v(p_0)^2}}
       =\operatorname{sgn}(u(p_0))u.
\]
Thus, for polygons, the differentiation argument in the proof of
Lemma~\ref{lem:strip-good-points} and the approximation of $p_0$ by
points $p_j$ in the proof of Proposition~\ref{prop:prescribed} are not
needed. The corner and selection arguments above are independent of
Theorem~\ref{thm:smooth}. The subsequent smooth approximation and
compact-superlevel argument are unchanged and still use that theorem
to obtain (HS1).

\subsection{Explicit starlike polygons with a simple eigenvalue}
\label{subsec:starlike}

The following family illustrates the simple-eigenvalue case with
nonconvex polygons. The notation of this subsection is independent of
the strip construction of Section~\ref{sec:lipschitz-bridge}.
Let $R=(-2,2)\times(-1,1)$ and let
$\vartheta:[-2,2]\to[0,1]$ be continuous and affine on each interval of
a finite partition, with $\operatorname{Lip}(\vartheta)\le1/2$. For
$0<\delta\le1/4$, set
\begin{equation}\label{eq:starlike-polygon-family}
 \omega(x)=1-\delta\vartheta(x),\qquad
 \Omega_{\delta,\vartheta}
 =\{(x,y):-2<x<2,\ |y|<\omega(x)\}.
\end{equation}
Since $\omega-x\omega'\ge1/2$ almost everywhere, the function
$t\mapsto\omega(tx)/t$ is nonincreasing on $(0,1]$, and hence
$\omega(tx)\ge t\omega(x)$ for
$0<t<1$. Thus $t\Omega_{\delta,\vartheta}\subset\Omega_{\delta,\vartheta}$,
which proves that these polygonal domains are starlike.
Their first positive Neumann eigenvalue is simple. We prove this
using the uniform estimate
\begin{equation}\label{eq:starlike-polygon-gap}
 \Lambda_{\Omega_{\delta,\vartheta},1}^{\mathrm N}
 \le\frac{\pi^2}{8}
 <\frac{\pi^2}{6}
 \le\Lambda_{\Omega_{\delta,\vartheta},2}^{\mathrm N}.
\end{equation}

Indeed, let $\Psi:R\to\Omega_{\delta,\vartheta}$,
$\Psi(x,y)=(x,\omega(x)y)$. For $f\in H^1(\Omega_{\delta,\vartheta})$
and $\tilde f=f\circ\Psi$ we have
$\int_{\Omega_{\delta,\vartheta}}f^2=\int_R\omega\tilde f^2$ and
$\int_{\Omega_{\delta,\vartheta}}|\nabla f|^2
=\int_R\mathsf M\nabla\tilde f\cdot\nabla\tilde f$, where
\[
 \mathsf M=\begin{pmatrix}
 \omega&-y\omega'\\ -y\omega'&(1+y^2\omega'^2)/\omega
 \end{pmatrix},\qquad
 \det\mathsf M=1,\quad \tfrac34\le\omega\le1,\quad |y\omega'|\le\tfrac18.
\]
The sum of the absolute values of the entries in each row of
$\mathsf M$ is at most $71/48$, so
$\lambda_{\max}(\mathsf M)\le71/48<3/2$.
Since $\det\mathsf M=1$,
$\lambda_{\min}(\mathsf M)=1/\lambda_{\max}(\mathsf M)\ge48/71>2/3$.
The min--max principle, together with
$\tfrac34\le\omega\le1$, yields, for every $k\ge0$,
\[
 \tfrac23\Lambda_{R,k}^{\mathrm N}
 \le\Lambda_{\Omega_{\delta,\vartheta},k}^{\mathrm N}
 \le\tfrac{71}{36}\Lambda_{R,k}^{\mathrm N}
 \le2\Lambda_{R,k}^{\mathrm N}.
\]
As $\Lambda_{R,1}^{\mathrm N}=\pi^2/16$ and
$\Lambda_{R,2}^{\mathrm N}=\Lambda_{R,3}^{\mathrm N}=\pi^2/4$, this proves
\eqref{eq:starlike-polygon-gap}.

\begin{figure}[htbp]
\centering
\begin{tikzpicture}[x=.92cm,y=.92cm,font=\small,
 polyoutline/.style={draw=blue!65!black,line width=.8pt,fill=blue!3},
 ray/.style={draw=orange!70!black,densely dashed,line width=.45pt}]
\begin{scope}
 \path[polyoutline] (-2,-1)--(0,-.75)--(2,-1)--(2,1)--(0,.75)--(-2,1)--cycle;
 \draw[ray] (0,0)--(0,.75) (0,0)--(0,-.75);
 \fill (0,0) circle (1.1pt);
 \node[right,inner sep=3pt] at (0,0) {$0$};
 \node at (0,1.4) {(a) Central indentations};
 \node at (0,-1.4) {$\vartheta=\vartheta_1$};
\end{scope}
\begin{scope}[xshift=4.6cm]
 \path[polyoutline] (-2,-1)--(0,-1)--(1,-.875)--(2,-1)
 --(2,1)--(1,.875)--(0,1)--(-2,1)--cycle;
 \draw[ray] (0,0)--(1,.875) (0,0)--(1,-.875);
 \fill (0,0) circle (1.1pt);
 \node[left,inner sep=3pt] at (0,0) {$0$};
 \node at (0,1.4) {(b) Off-center indentations};
 \node at (0,-1.4) {$\vartheta=\vartheta_2$};
\end{scope}
\begin{scope}[xshift=9.2cm]
 \path[polyoutline] (-2,-1)--(-1,-.875)--(0,-1)--(1,-.875)--(2,-1)
 --(2,1)--(1,.875)--(0,1)--(-1,.875)--(-2,1)--cycle;
 \draw[ray] (0,0)--(1,.875) (0,0)--(1,-.875) (0,0)--(-1,.875) (0,0)--(-1,-.875);
 \fill (0,0) circle (1.1pt);
 \node[above,inner sep=3pt] at (0,0) {$0$};
 \node at (0,1.4) {(c) Four indentations};
 \node at (0,-1.4) {$\vartheta=\vartheta_3$};
\end{scope}
\end{tikzpicture}
\caption{Nonconvex starlike polygons from
\eqref{eq:starlike-polygon-family} with $\delta=1/4$, drawn to scale.
The closure of each polygon lies in the rectangle
$[-2,2]\times[-1,1]$ and contains its four corners.
The polygons are starlike with respect to
the marked origin. In every panel, $\dim E_1=1$ by
\eqref{eq:starlike-polygon-gap}. The dashed segments join the origin to
the reentrant vertices and lie in the closure of the polygon.}
\label{fig:starlike-simple}
\end{figure}
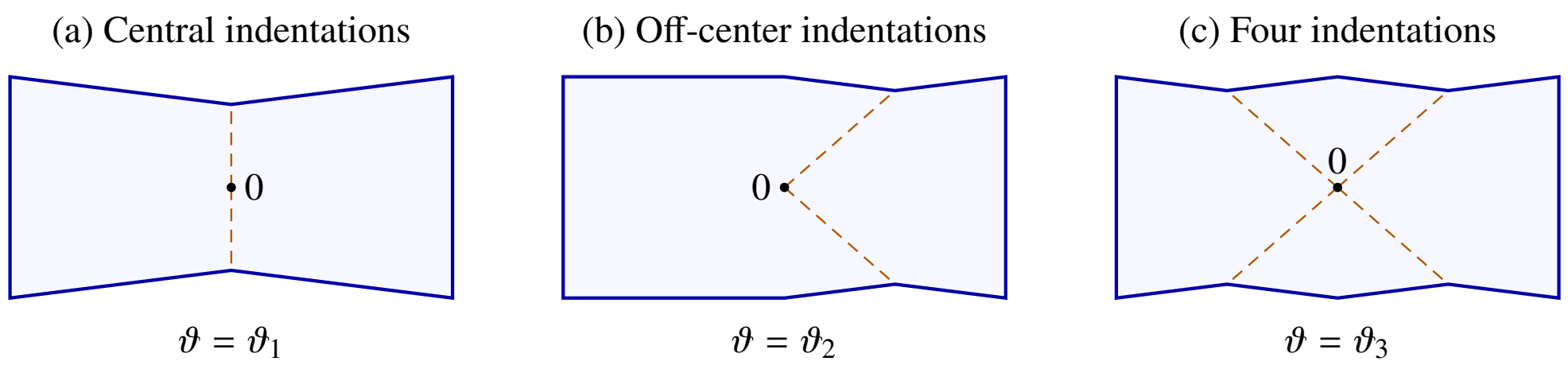

Figure~\ref{fig:starlike-simple} shows the nonconvex choices
\begin{equation}\label{eq:starlike-polygon-profiles}
 \vartheta_1(x)=1-\frac{|x|}{2},\qquad
 \vartheta_2(x)=\frac{(1-|x-1|)_+}{2},\qquad
 \vartheta_3(x)=\frac{\dist(x,\{-2,0,2\})}{2},
\end{equation}
where $s_+=\max\{s,0\}$. Each satisfies the stated bounds.
Moreover, $\omega(\pm2)=1$ and $\omega<1$ somewhere, so $\omega$ is
not concave and the corresponding polygon is nonconvex.
The construction and the spectral-gap estimate
\eqref{eq:starlike-polygon-gap} are independent of
Theorems~\ref{thm:lipschitz} and~\ref{thm:smooth}.
Their (HS1) property follows from Theorem~\ref{thm:lipschitz}, whose
proof uses the smooth input.

These three families are not lip domains in the sense of
Atar and Burdzy \cite{AB}, even after rigid motions. The sides of a lip
domain have unoriented directions in a closed angular interval of
width $\pi/2$. The two graphs that bound a lip domain meet at its
leftmost and rightmost points, so in the defining coordinates every
side has slope in $[-1,1]$. The directions of a vertical side and of sides of slopes
$\delta/2$ and $-\delta/2$ do not lie in any such interval.

\FloatBarrier

\begin{remark}
Using OpenAI's Theorem~\ref{thm:smooth} as an input,
K.~Burdzy generated, with the assistance of ChatGPT, an argument
establishing the HS1 property for compact simply connected
surfaces with smooth boundary and a smooth flat metric.
The argument has not been verified by a human. This extension concerns
the geometry of the underlying surface, whereas the present
paper concerns boundary regularity. 
\end{remark}

\section*{Acknowledgments}
The author gratefully acknowledges Krzysztof Burdzy for many valuable discussions on hot spots over the years, particularly in recent months.  The author thanks Greg Buzzard for many conversations on the topic of this paper and for obtaining a review of the manuscript from Claude (Anthropic), which offered interesting insights and raised further questions. He also thanks Greg and Chris Janjigian for many conversations on the use of AI in the mathematical sciences and STEM fields more broadly.

\section*{AI Disclosure}
This manuscript was prepared with assistance from ChatGPT (OpenAI)
in exposition, calculations, and proof checking. The approximation
approach and constructions are due to the author. The smooth-domain
theorem is due to OpenAI and is used here as a cited input. For an extended  expostion of the OpenAI proof, see also \cite{BanuelosSmooth2026}.
The classical ingredients and the critical-set argument of Judge and
Mondal are identified where they are used. The author takes
responsibility for the mathematical content of the paper.

\end{document}